\documentclass[11pt, a4paper]{elsarticle}
\usepackage[utf8]{inputenc}
\makeatletter
\def\ps@pprintTitle{%
   \let\@oddhead\@empty
   \let\@evenhead\@empty
   \let\@oddfoot\@empty
   \let\@evenfoot\@oddfoot
}
\makeatother

\usepackage{setspace}
\usepackage{amsmath,hhline,dcolumn,layout,booktabs, pdflscape}
\usepackage{algorithm}
\usepackage{algpseudocode}
\usepackage{floatflt}
\usepackage{graphicx}
\usepackage{subcaption}
\usepackage{longtable,verbatim}
\usepackage{enumitem,  rotfloat, supertabular}
\usepackage{lscape}
\usepackage{rotating}
\usepackage{booktabs}
\usepackage{amsfonts}
\usepackage{multirow}
\usepackage{chapterbib}
\usepackage{hyphenat}
\usepackage{multirow}
\usepackage{url}
\usepackage{xcolor}
\usepackage{mathtools}
\usepackage{amssymb}
\usepackage{verbatim}

\newtheorem{corollary}{Corollary}[section]
\newtheorem{proposition}{Proposition}[section]
\allowdisplaybreaks

\title{
{Modeling and Optimizing Resource Allocation Decisions  through Multi-model Markov Decision Processes with Capacity Constraints}}

\author[koc]{Onur Demiray\corref{cor1}}
\author[koc2]{Evrim Didem Güneş}
\author[koc]{Lerzan Örmeci}
\address[koc]{Department  of  Industrial  Engineering,  Koç  University,  34450 Istanbul, Turkey}
\address[koc2]{College  of  Administrative  Sciences  and  Economics,  Koç  University, 34450 Istanbul, Turkey}

\cortext[cor1]{Onur Demiray, odemiray18@ku.edu.tr, Department  of  Industrial  Engineering,  Koç  University,  34450 Istanbul, Turkey}

\begin{document}

\begin{abstract}
\small

This paper proposes a new formulation for the dynamic resource allocation problem, which converts the traditional MDP model with known parameters and no capacity constraints to a new model with uncertain parameters and a resource capacity constraint.
Our motivating example comes from a medical resource allocation problem: patients with multiple chronic diseases can be provided either normal or special care, where the capacity of special care is limited due to financial or human resources. In such systems, 
it is difficult, if not impossible,  to generate good estimates for the evolution of health for each patient.
We formulate the problem as a two-stage stochastic integer program. However, it becomes easily intractable in larger instances of the problem for which we propose and test a parallel approximate dynamic programming algorithm. We show that commercial solvers are not capable of solving the problem instances with a large number of scenarios. Nevertheless, the proposed algorithm provides a solution in seconds even for very large problem instances. In our computational experiments, it finds the optimal solution for $42.86\%$ of the instances. On aggregate, it achieves $0.073\%$ mean gap value. Finally, we estimate the value of our contribution for different realizations of the parameters. Our findings show that there is a significant amount of additional utility contributed by our model.

\end{abstract}


\maketitle
\onehalfspacing

\section{Introduction}
\label{section:introduction}

Markov Decision Processes (MDPs) are successfully used to find optimal policies in sequential decision making problems under uncertainty. The application areas of MDPs vary from inventory management, finance, robotics, telecommunication to humanitarian logistics \cite{altman1999constrained, boucherie2017markov, merakli2019risk,puterman2014markov}. It is also extensively employed in medical decision making literature. Practices in this literature range from determining the initiation time of a drug to scheduling patients for a treatment \cite{alagoz2004optimal, denton2009optimizing, denton2018optimization, shechter2008optimal}. We refer to \cite{alagoz2010markov} 
for a comprehensive review on the medical practices of MDPs. However, they have certain limitations, as they assume that both transition probability matrices for all actions and the rewards incurred for all state-action pairs are known. In practice, these parameters are generally estimated by maximum likelihood estimators derived using observational data \cite{zhang2017robust}. However, having the right data set and deriving reliable estimators are notoriously hard and it is prone to errors which may cause suboptimal solutions \cite{denton2018optimization, grand2020robust, mannor2016robust, mannor2007bias, merakli2019risk, zhang2017robust, zhang2017probabilistic}. Hence, in recent years, there are studies on  MDP models with uncertain parameters, such as \cite{steimle2018multi}. Another important feature that traditional MDP models generally ignore is the capacity constraints, which has also received increasing attention in the recent years \cite{ayvaci2012effect, cevik2018analysis, deo2013improving}. Our work combines these two emerging areas, as we consider a system with uncertain transition probability matrices and rewards, which operate under a resource constraint. 

Robust optimization is extensively employed to deal with parameter uncertainty in MDPs \cite{grand2020robust, merakli2019risk, zhang2017robust}. Robust optimization is one of the widely used techniques for optimization under uncertainty where no probability distribution is assumed for uncertain data. 
It is based on  the worst-case realization of the parameters from a set of alternatives what is so-called uncertainty set \cite{ben2009robust,ben2002robust,bertsimas2011theory, gorissen2015practical}. 
A number of studies  approached the parameter uncertainty in MDP problem employing robust optimization with a polyhedral uncertainty set due to its properties allowing tractable solution algorithms \cite{givan2000bounded,satia1973markovian,tewari2007bounded,white1994markov}, whereas another group 
modeled the problem with the rectangularity assumption which ensures the independence of the rows in transition probability matrices \cite{iyengar2005robust, nilim2004robustness, sinha2016policy}. Wiesemann et al. \cite{wiesemann2013robust} studied robust MDPs by relaxing this assumption.

Parameter uncertainty in MDPs is also examined with stochastic programming where different types of Bayesian approaches are adopted \cite{merakli2019risk}. In this context, Steimle et al. \cite{steimle2018multi} propose a multi-model MDP (MMDP) approach for a finite-horizon MDP, where transition probabilities and rewards can come from different  models. They further claim that having a single policy optimized with respect to all models can be a good idea to capture variations within MDP parameters. An MMDP model can be seen as a two-stage stochastic programming formulation where different models correspond to different scenarios. Meraklı $\&$ Küçükyavuz \cite{merakli2019risk} provided a stochastic programming model for an infinite-horizon MDP model under a risk-averse setting, and applied it to a humanitarian inventory management problem.

Studies focusing on parameter uncertainty in MDPs in the literature address only the unconstrained version of the problem. Nevertheless, some healthcare applications employing MDPs necessitate taking capacity constraints into account due to existance of resource scarcity \cite{ayvaci2012effect, cevik2018analysis, deo2013improving}. Deo et al. \cite{deo2013improving} proposed a constrained MDP model scheduling patients for a community-based chronic care delivery problem. Ayvaci et al. \cite{ayvaci2012effect} and Cevik et al. \cite{cevik2018analysis} developed constrained MDP and partially observable MDP (POMDP) models of diagnostic and screening decisions for breast cancer. However, none of these studies incorporated parameter uncertainty which may cause suboptimal solutions for such medical resource allocation problems.

With a motivation from medical resource allocation problems, we, in this paper, introduce a capacity constrained, finite-horizon, and finite-state MDP model with uncertain transition probabilities and rewards. In addition to the finite number of non-absorbing states, we consider only one absorbing state to reflect all conditions that stop the decision process. Under the healthcare setting, this absorbing state corresponds to all the terminating conditions such as death, stroke, or hearth failure. The another assumption we make in our model is that the action space consists of two actions, which is consistent considering the nature of the resource allocation problems (allocate versus don't allocate or accept versus reject). Though we restrict the model to have two actions in each decision epoch, the formulation can be extended easily to the cases with more complicated action spaces. Moreover, the model can represent any capacity constrained problem using MDPs, i.e., its scope is larger than  medical resource allocation problems. For example, 
our approach can be applied to the dynamic product portfolio management problem (Seifert et al. \cite{seifert2016dynamic}) under a financial constraint, where the probabilities associated with product life cycle transitions are uncertain. The case of a varying level of working capital over different business cycles can also be addressed.  Thus the proposed model can be used under different domains spanning from finance, inventory management to marketing.

We adopt a multi-model MDP approach to solve our model. In other words, we provide an extensive-form formulation that is an approximation of the underlying two-stage stochastic integer program with the help of multiple MDP representations, i.e., scenarios. In this respect, our model can be considered as an extension of the model developed by Steimle et al. \cite{steimle2018multi} with capacity constraints. They formulated their models based on the primal linear program (LP) developed for finding optimal policies for an unconstrained MDP. However, we build our formulation based on the dual LP which is the version allowing us to incorporate additional constraints to the MDP framework. To the best of our knowledge, this is the first study combining capacity constraints and parameter uncertainties in MDPs.

The extensive-form formulation, i.e., the deterministic equivalent of the two-stage stochastic programming formulation, can easily become intractable because of capacity constraints, scenarios, and decision epochs. With this in mind, we first prove results leveraging the problem structure, and then develop a parallel approximate dynamic programming algorithm based on these results. This algorithm makes us capable of solving the model efficiently even if we have a large number of scenarios and decision epochs. Lastly, the proposed model and the corresponding solution approaches are applied to a chronic care delivery problem, which also forms the basis for the computational experiments. Through computational studies, we compare different solution approaches. We also measure the value of perfect information and stochastic solution introduced by taking transition probability and reward uncertainties into account. It is worth to mention that all experiments are conducted for a fixed state space which provides a reasonable number of states for our analysis. Thus, our algorithm reflects to the challenges for increasing number scenarios and decision epochs.

Our main contributions can be summarized as follows:
\begin{itemize}
    \item We introduced a new model incorporating both parameter uncertainty and capacity constraints  in MDPs. To the best of our knowledge, this is the first study considering these two dimensions together.
    \item We prove problem specific structural results, and propose a novel network representation, which together enable a parallel approximate dynamic programming algorithm to solve larger instances of the problem. 
    \item We conduct extensive computational experiments to measure the value of embracing parameter uncertainty and to compare the solution approaches using the chronic care delivery problem setting as an example.
\end{itemize}

The remainder of this paper is organized as follows. The model is formally described and mathematically modeled in Section \ref{section:definition}. Section \ref{section:approach} explains the proposed solution approach. A comprehensive computational study is provided in Section \ref{section:experiment} after introducing the chronic care delivery problem in Section \ref{section:case}. Finally, Section \ref{section:conclusion} concludes the paper and draws directions for future research.

\section{Model Description and Formulation}
\label{section:definition}

We consider a capacity allocation problem in a system with  $N$ homogeneous individuals who can be in different states throughout the planning horizon. 
The actions are taken at the population level to respect  a capacity constraint, whereas the individuals' states evolve depending on their initial states and on the actions. We first describe the process for an individual, and then relate it with the process at the  population level.

The state of an individual evolves according to an uncertain transition probability matrix. In addition, the rewards collected are also uncertain. Hence, her evolution  can be represented by a Markov decision process with uncertain parameters, which corresponds to  a multi-model Markov decision process (MMDP) with the corresponding tuple $(\mathcal{T}, \mathcal{S}, \mathcal{A}, \Omega, \Lambda)$ (see e.g., \cite{steimle2018multi}). Here,  $\tilde{\mathcal{T}} \triangleq \{1, \dotsc, T-1\}$ is the set of decision epochs and $\mathcal{T}=\tilde{\mathcal{T}} \cup \{T\}$ is the set of all periods; an individual can be in state $s\in\mathcal{S}$ in time period $t$,  where $\mathcal{S}$ is the set of all states, $\tilde{\mathcal{S}}$ is the set  of non-absorbing states, and  $D$ is the  single  absorbing state, so that $\mathcal{S}=\tilde{\mathcal{S}} \cup \{D\}$; $\mathcal{A}=\{0,1\}$ is the set of action space, where $a=1$ corresponds to special service with a limited capacity, and $a=0$ to regular service with unlimited capacity ; $\Omega$ is the set of models or scenarios that specify different transition probability distributions and rewards; and ${\Lambda}=(\lambda_{\omega})_{\omega\in\Omega}$, where $\lambda_{\omega}>0$ is the probability that scenario $\omega$ is in effect and $\sum_{\omega\in \Omega}\lambda_{\omega}=1$. 

Scenario $\omega$ specifies two transition probability matrices; one among the non-absorbing states, $\mathcal{P}^\omega \triangleq \{P_{iaj}^\omega\}_{(i,j)\in \tilde{\mathcal{S}} \times  \tilde{\mathcal{S}}, a \in \{0,1\}}$, and the other from non-absorbing states to the absorbing state $D$, $\mathcal{Q}^\omega \triangleq \{Q_{ia}^\omega\}_{i \in \tilde{\mathcal{S}}, a \in \{0,1\}}$. 
In words, when action $a\in\mathcal{A}$ is taken in state $i\in\mathcal{\tilde{S}}$ 
under scenario $\omega \in \Omega$, 
$P_{iaj}^\omega$ is the probability that an individual in state $i$ moves to state $j \in \tilde{\mathcal{S}}$, whereas $Q_{ia}^\omega$ is the probability that the individual reaches the absorbing state.
For each scenario $\omega \in \Omega$, $\mathcal{P}
^\omega$ and $\mathcal{Q}^\omega$ satisfy the following inequalities:
\begin{subequations}
\begin{align}
    P_{iaj}^\omega \geq 0  && \forall (i,j)\in \tilde{\mathcal{S}} \times \tilde{\mathcal{S}}, a \in \mathcal{A}, \omega \in \Omega \label{validstochastic} \\
    Q_{ia}^\omega \geq 0 && \forall i \in \tilde{\mathcal{S}}, a \in \mathcal{A}, \omega \in \Omega \label{validstochastic1.5}\\
     \sum_{j \in \tilde{\mathcal{S}}} P_{iaj}^\omega + Q_{ia}^\omega = 1 && \forall i \in \tilde{\mathcal{S}}, a \in \mathcal{A}, \omega \in \Omega \label{validstochastic2}
\end{align}
\end{subequations}
For all non-absorbing states $i \in \tilde{\mathcal{S}}$, the reward gained 
by action $a$  under scenario $\omega$ is denoted by $r_{ia}^\omega$, whereas the reward  
of being in state $i$ at the final stage $T$ is given by $R_i^\omega$.
We further assume that the reward of visiting the absorbing state, $R^D$, is fixed for all scenarios. 

In general, the objective of MMDPs is to find a strategy over the entire planning horizon, denoted $\Pi$, that performs well with respect to different scenarios accounted in MMDP. One of the widely-used techniques to solve such problems is to employ a two-stage stochastic program, where the first stage determines the strategy $\Pi$ and   the second stage evaluates the performance of $\Pi$. 
The discrete set of scenarios called \textit{solution sample} allows developing a mixed integer program that solves for both stages at the same time, which is also called {\it extensive form of MMDP}. Our solution approach also uses this technique; however it has to  incorporate the additional capacity constraint at the population level, as we describe below.

The initial population consists of $N$ individuals, where $n_i$ of them are in state $i$ so that $N=\sum_{i \in \tilde{\mathcal{S}}}n_i$. Consequently, the initial distribution of states is specified by $\theta_i = \frac{n_i}{N}$. 
The state of the system in period $t$ is given by $\mathbf{n}(t)=(n_i(t))_{i\in\tilde{\mathcal{S}}}$, where $N-\sum_{i \in \tilde{\mathcal{S}}}n_i(t)$ individuals are in the absorbing state at time $t$.
Individuals receive either regular or special service, where the special service needs additional resources. Hence, the problem is a capacity allocation problem, where scarce resources are rationed among the individuals in different states in each decision period to maximize the total expected reward over $T$ periods.
We represent the  scarcity of the resources by a constraint that limits the expected number of individuals who receive special service by $C_t$ at decision epoch $t \in \tilde{\mathcal{T}}$. This constraint requires tracking the evolution of  all individuals through the states as well as of the actions over time.

 The system manager decides on the group(s) of individuals who will receive special service in each period, i.e., she chooses a subset of $\tilde{\mathcal{S}}$ who are entitled to use the additional resources under this constraint.
Hence, the action space in all states is given by $\mathcal{D}=(0,1)^{|\tilde{\mathcal{S}}|}$.
We set $\pi^t_i=1$ if the resource is used for the individuals in state $i$ at decision epoch $t$, and $\pi^t_i=0$ otherwise, where $i \in \tilde{\mathcal{S}}$. 
The policy in decision period $t$ is then given by  $\pi^t=(\pi^t_i)_{i\in\tilde{\mathcal{S}}}$, and the strategy over the planning horizon by
$\Pi \triangleq (\pi^t)_{t \in \tilde{\mathcal{T}}}$.
This study aims to find a deterministic strategy $\Pi$ that maximizes the total expected reward over $T$ periods, while satisfying the capacity constraints in each decision period under all  scenarios.
Hence, the optimal strategy $\Pi$, which 
is independent of the scenarios, 
determines the groups that will receive special service in all periods $t\in\tilde{\mathcal{T}}$.

To account for the expected rewards and to respect the capacity constraints, we represent the evolution of the system in the model by defining  the so-called \textit{occupancy measures} for each scenario $\omega$:
\begin{longtable*}{ p{.03\textwidth}  p{.97\textwidth}}
$\mathcal{X}_{ia}^{\omega t}$ & probability of being in state $i$ under action $a$ at decision epoch $t$\\
$\mathcal{Z}^{\omega t}$ & probability of visiting the absorbing state at period $t$\\
$\mathcal{Y}_{i}^\omega$ & probability of finalizing the decision process in non-absorbing state $i$
\end{longtable*}

For a given strategy $\Pi$, the total expected reward is computed as follows:
\begin{equation}
  \mathcal{U} (\Pi)= \underbrace{N \underbrace{\Bigg [\sum_{\omega \in \Omega} \lambda_\omega \Big(
    \underbrace{\overbrace{\sum_{t \in \mathcal{T} \setminus \{1\}}\mathcal{Z}^{\omega t}R^D}^{\text{component 1}}
    + \overbrace{\sum_{t \in \tilde{\mathcal{T}}}\sum_{i \in \tilde{\mathcal{S}}}  \mathcal{X}^{\omega t}_{i\pi^t_i}r_{i\pi^t_i}^\omega}^{\text{component 2}} +\overbrace{ \sum_{i \in \tilde{\mathcal{S}}} \mathcal{Y}_{i}^{\omega}R_i^\omega}^{\text{component 3}}}_{\text{expected reward of an individual under scenario $\omega \in \Omega$}} 
    \Big) \Bigg].}_{\text{expected reward of an individual}}}_{\text{total expected reward of all individuals}} \label{totalExpectedReward}
\end{equation}
The quantity in the brackets corresponds to the expected reward of an individual under all scenarios, which gives the total expected reward of the system when multiplied by the total population $N$. The expected reward of an individual under scenario $\omega \in \Omega$ has three components:
Component 1 computes the expected reward over all periods due to visits to the absorbing state. Component 2 accounts for the expected reward of an individual in the non-absorbing states throughout all stages $t\in\tilde{\mathcal{T}}$. Finally, component 3 determines the expected reward in the final period, $T$.

At this point, we point out certain characteristics of the problem:
(1) Actions $\pi_i^t$ can be either 0 or 1, so that either all individuals in state $i$  use the additional resources or all are excluded from the special service. 
We can easily relax this condition  to allow partial coverage,  
as demonstrated in some of our numerical experiments in Section 5.
(2) Policy $\pi^t$ changes with respect to time $t$, but  not with respect to the state of the system, 
$\mathbf{n}(t)$. More explicitly, we determine the strategy $\Pi$ for the entire planning horizon at time 0 and 
commit to it till the end of period $T$. This allows us to account for the future effects of our current actions 
in expectation. Note that in practice, the realized number of individuals to receive the special service may exceed the capacity. 
However in this paper we aim to support strategic-level decisions, whereas such violations which may 
be encountered and dealt with at the operational level are ignored.
Solving this model on a rolling-time horizon basis by  replacing the initial population with the current number of individuals in each period will account for the dynamic change of the system state and ensure that the capacity constraint is never violated.


The mixed integer programming (MIP) formulation demonstrated through (\ref{ef_obj})-(\ref{RNEF_domain4}) finds $\Pi^*$ which maximizes the total expected reward of the population without violating the capacity constraints for the MMDP model. For the remainder of the paper, we call it (MIP-MMDP). It is the deterministic equivalent of the underlying two-stage stochastic integer programming formulation. In this formulation, $\Pi=\{\pi_i^t\}_{i \in \tilde{\mathcal{S}}, t \in \tilde{\mathcal{T}}}$ form the first-stage decision variables that need to be made before nature reveals the uncertainty. Hence, they are independent of scenarios. The second-stage problem, on the other hand, determines the values of the occupancy measures after values of the uncertain parameters are realized by the scenarios. As mentioned above, these two stages can be combined in a single mixed-integer program when solving for MMDPs.

\begin{align}
            & {\textbf{max}}
            & &  \mathcal{U}(\Pi) \label{ef_obj}\\
            & \textbf{st;}
            & & \sum_{a \in \mathcal{A}}\mathcal{X}^{\omega,1}_{ia}=\theta_i & \forall i \in \tilde{\mathcal{S}}, \omega \in \Omega \label{RNEF_flow1}\\
            &&& \sum_{i \in \tilde{\mathcal{S}}} \sum_{a \in \mathcal{A}} \mathcal{X}^{\omega,t-1}_{ia}P^{\omega}_{iaj}=\sum_{a \in \mathcal{A}} \mathcal{X}_{ja}^{\omega, t} & \forall t \in \tilde{\mathcal{T}}\setminus\{1\},  j \in \tilde{\mathcal{S}}, \omega \in \Omega \label{RNEF_flow2} \\
             &&& \sum_{i \in \tilde{\mathcal{S}}} \sum_{a \in \mathcal{A}} \mathcal{X}^{ \omega,T-1}_{ia}P^{\omega}_{iaj}=\mathcal{Y}^{\omega}_j & \forall j \in \tilde{\mathcal{S}}, \omega \in \Omega \label{RNEF_flow3} \\
             &&& \sum_{i \in \tilde{\mathcal{S}}}\sum_{a \in \mathcal{A}}\mathcal{X}^{\omega,t-1}_{ia}Q^{\omega}_{ia}=\Delta\mathcal{Z}^{\omega, t} &  \forall t \in \mathcal{T}\setminus\{1\}, \omega \in \Omega \label{RNEF_flow4}\\
             &&& N\sum_{i\in \tilde{\mathcal{S}}}\mathcal{X}_{i1}^{\omega, t}\leq C_t & \forall t \in \tilde{\mathcal{T}}, \omega \in \Omega \label{RNEF_cap}\\
             &&& \mathcal{X}_{i1}^{\omega, t} \leq \pi_{i}^{t} &\forall i \in \tilde{\mathcal{S}}, t \in \tilde{\mathcal{T}}, \omega \in \Omega  \label{RNEF_deterministicPol1}\\
            &&& \mathcal{X}_{i0}^{\omega, t} \leq (1-\pi^{t}_{i}) & \forall i \in \tilde{\mathcal{S}}, t \in \tilde{\mathcal{T}}, \omega \in \Omega  \label{RNEF_deterministicPol2}\\
            &&& \pi_{i}^{t} \in \{0,1\} & \forall t \in  \tilde{\mathcal{T}}, i \in  \tilde{\mathcal{S}}
            \label{RNEF_domain4}\\
            &&& \mathcal{X}_{ia}^{\omega, t} \geq 0 &\forall (t,i, a, \omega) \in ( \tilde{\mathcal{T}},  \tilde{\mathcal{S}}, \mathcal{A}, \Omega) \label{RNEF_domain1}\\
            &&& \mathcal{Z}^{\omega,t} \geq 0 & \forall t \in \mathcal{T}\setminus\{1\}, \omega \in \Omega \label{RNEF_domain2}\\
            &&& \mathcal{Y}^{\omega}_{i} \geq 0 &  \forall i \in  \tilde{\mathcal{S}}, \omega \in \Omega \label{RNEF_domain3}
\end{align}
where  $\Delta\mathcal{Z}^{\omega, t}=\mathcal{Z}^{\omega, t}-\mathcal{Z}^{\omega, t-1} \geq 0$,  and $\mathcal{Z}^{\omega, 1}=0$, for all $ \omega \in \Omega$.

In (MIP-MMDP), the objective function is given by (\ref{ef_obj}), which maximizes the total expected reward obtained by all individuals, as defined by Equation (\ref{totalExpectedReward}). Equations in (\ref{RNEF_flow1})-(\ref{RNEF_flow4}) are called \textit{forward equations}, and they recursively compute the values of the occupancy measures. From a different viewpoint, they balance the probability flows in the model. Furthermore, (\ref{RNEF_cap}) ensures that capacity constraints are satisfied for each scenario. Particularly, the expected number of individuals entitled to special service, i.e., those in states with action $a=1$, cannot exceed the allowed capacity $C_t$ at decision epoch $t \in \tilde{\mathcal{T}}$ for each scenario. Constraints (\ref{RNEF_deterministicPol1}) and (\ref{RNEF_deterministicPol2}) 
establish the link between the first stage decisions $\Pi$ and the second stage variables $\mathcal{X}$. Specifically, if the strategy $\Pi$ does not take action $a$ for the group of individuals in state $i$ at time $t$, then the corresponding occupancy measures has to be equal to 0. Otherwise, they can take values up to 1. 
Constraints (\ref{RNEF_domain4}) ensure that all individuals in state $i$ receive either regular or special service. 

We now  provide  two  valid  inequalities for  the  (MIP-MMDP) formulation. Their proofs  are presented in \ref{appendix:proof_valid}. Proposition~\ref{VI:1} is a simple consequence of transition probability distributions,
while Proposition \ref{VI:2} is achieved by aggregating the capacity constraints over all decision epochs.

\begin{proposition}
\label{VI:1}
Hyperplanes in (\ref{eq_VI_1}) are valid inequalities for (MIP-MMDP).
\begin{align}\label{eq_VI_1}
1 = 
        \begin{dcases}
            \sum_{i \in \tilde{\mathcal{S}}} \sum_{a \in \mathcal{A}} \mathcal{X}_{ia}^{\omega, t}, & \text{if } t=1\\
             \sum_{i \in \tilde{\mathcal{S}}} \sum_{a \in \mathcal{A}} \mathcal{X}_{ia}^{\omega, t} + \mathcal{Z}^{\omega, t}, & \text{if } t \in \{2 ,\dotsc, T-1\} \\
             \sum_{i \in \tilde{\mathcal{S}}}\mathcal{Y}_i^\omega + \mathcal{Z}^{\omega, t}, & \text{if } t = T
        \end{dcases}\
        \quad \forall \omega \in \Omega
\end{align}
\end{proposition}

\begin{proposition}
\label{VI:2}
$\sum\limits_{t \in \tilde{\mathcal{T}}} \Big [\sum\limits_{i \in \tilde{\mathcal{S}}}\mathcal{X}^{\omega, t}_{i,0} + \mathcal{Z}^{\omega, t}\Big] \geq T - \frac{\sum\limits_{t \in \tilde{\mathcal{T}}}C_t + n}{n}$ defines a valid inequality for each scenario $\omega \in \Omega$ for (MIP-MMDP).
\end{proposition}

\section{Solution Approach}
\label{section:approach}

Global optima may be obtained in reasonable times for (MIP-MMDP) with small and medium-size problem instances by using commercial solvers. However, reaching the optimal solution(s) for larger instances of the problem is typically beyond the realm of possibility. Nevertheless, the applicability of the model for real settings requires a huge number of parameters such as a great number of scenarios in order to satisfy the convergence of the extensive-form formulation or many decision epochs. Therefore, we propose a parallel approximate dynamic programming (PADP) algorithm providing optimal solutions or very close solutions to global optima as we will show through computational experiments. To that end, we first prove structural results by which we are able to represent our problem in a specific graph structure. Then, we provide the details of the algorithm.

\subsection{Structural Properties}
\label{subsection:structure}

\begin{proposition}
\label{prop:uniqueness}
Given a strategy $\Hat{\Pi}$, the values of the occupancy measures under $\Hat{\Pi}$, denoted by $(\mathcal{X},\mathcal{Y}, \mathcal{Z})^{\Hat{\Pi}}$, are unique.
\end{proposition}

Proposition \ref{prop:uniqueness} is a consequence of the fact that when a strategy is fixed, then the problem is reduced to a Markov Reward Process (MRP). That is, $(\mathcal{X},\mathcal{Y}, \mathcal{Z})^{\Hat{\Pi}}$ is a result of the strategy evaluation algorithm illustrated in Algorithm \ref{Algorithm:MRP} (see \ref{appendix:mrp}) which utilizes the forward equations in (\ref{RNEF_flow1})-(\ref{RNEF_flow4}) to find the corresponding occupancy measures. Algorithm \ref{Algorithm:MRP} is a polynomial time algorithm in which computations for each scenario are tackled with parallel processing. 

\begin{corollary}\label{cor:feasibility}
The feasibility of a given strategy $\Hat{\Pi}$ can be determined by Algorithm \ref{Algorithm:MRP}.
\end{corollary}

\begin{corollary}\label{cor:forward}
Consider two strategies $\hat{\Pi}$ and $\bar{\Pi}$ with the following property: There exists $t^* \in \tilde{\mathcal{T}}\setminus\{T-1\}$ such that $\hat{\pi}_i^t = \bar{\pi}_i^t$ for all $i \in \tilde{\mathcal{S}}, t\leq t^*$, and $\hat{\pi}_{t^*+1}^i \neq \bar{\pi}_{t^*+1}^i$ for some $i \in \tilde{\mathcal{S}}$. Even though $\hat{\Pi} \neq \bar{\Pi}$, the occupancy measures until $t^*$ under $\hat{\Pi}$ and $\bar{\Pi}$ equal each other.
\end{corollary}

Corollary \ref{cor:forward} stems from the fact that Algorithm \ref{Algorithm:MRP} computes the occupancy measures in a forward direction. Thus, any change in policy $\pi^t$ only impacts the occupancy measures related to the stage $t$ and afterwards. This idea is the cornerstone of the proposed algorithm, as it enables us to express the problem in a network structure.

\subsection{The Parallel Approximate Dynamic Programming Algorithm}
\label{subsection:adp}

The network can be conceived as a 2-Dimensional topology in which the horizontal axes is constituted by the stages in $\tilde{\mathcal{T}}$, whereas all policy combinations, $\Upsilon \triangleq \{1,\dotsc, \upsilon,\dotsc,2^{|\tilde{\mathcal{S}}|}\}$, form the vertical axes. Then, we generate a node for each element in $\tilde{\mathcal{T}} \times \Upsilon$. A node associated with stage $t \in \tilde{\mathcal{T}}$ and combination $\upsilon \in \Upsilon$ is represented by $\pi_t(\upsilon)$. Visiting the node $\pi_t(\upsilon)$ means that the policy encoded by the combination $\upsilon\in\Upsilon$ is accepted for stage $t \in \tilde{\mathcal{T}}$ as part of the strategy. 

\begin{figure}[t]
\centering
\caption{Network Representation}
\includegraphics[width=\linewidth]{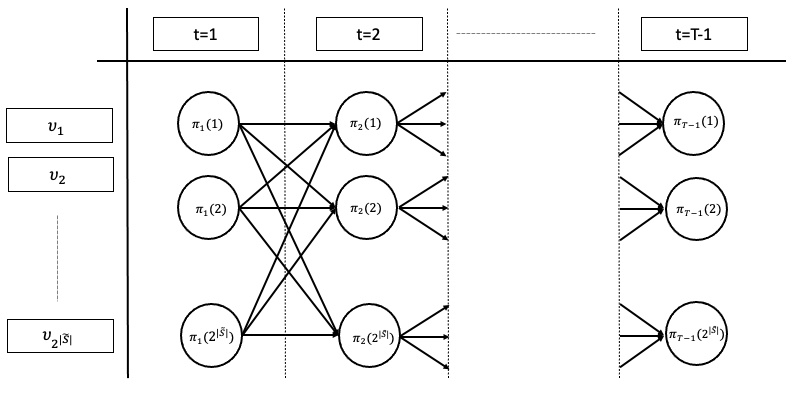}
\label{fig:exampleNet}
\end{figure}

$\phi_t \triangleq \{\pi_t(\upsilon)\}_{\upsilon\in\Upsilon}$ shows the set of nodes associated with stage $t \in \tilde{\mathcal{T}}$. Then, $\Phi \triangleq \{\phi_t\}_{t \in \tilde{\mathcal{T}}}$ becomes the set of nodes in the graph. Arcs take place only between the nodes of two consecutive stages. An arc from $\pi_t(\upsilon_1) \in \phi_t$ to $\pi_{t+1}(\upsilon_2) \in \phi_{t+1}$ means that the strategy ensures the policy represented by combination $\upsilon_2 \in \Upsilon$ for stage $t+1 \in \tilde{\mathcal{T}}$ while the policy implied by $\upsilon_1 \in \Upsilon$ is accepted for stage $t \in \tilde{\mathcal{T}}$. Then, $\gamma_t \triangleq \{(\pi_t(\upsilon_1),\pi_{t+1}(\upsilon_2)):(\upsilon_1,\upsilon_2) \in \Upsilon \times \Upsilon\}$ shows the set of arcs outflowing from the nodes in $\phi_t$, $t \in \tilde{\mathcal{T}}\setminus\{T-1\}$, and $\Gamma \triangleq \{\gamma_t\}_{t \in \tilde{\mathcal{T}}\setminus\{T-1\}}$ becomes the set of arcs. Then, our graph can be defined as $\mathcal{G}=(\Phi, \Gamma)$. A \textit{strategy path} denoted by $\sigma$ is nothing but an ordered list of the nodes in $\Phi$ so that there is exactly one node from $\phi_t, \forall t \in \tilde{\mathcal{T}}$. 

\begin{figure}[ht]
\centering
\caption{A Simple Illustrative Example}
\includegraphics[width=\linewidth]{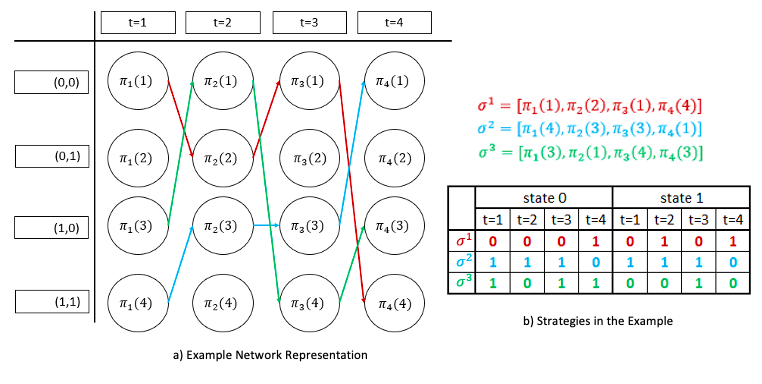}
\label{fig:exampleNetStrategy}
\end{figure}

Figure \ref{fig:exampleNetStrategy} provides an illustrative example with 2 non-absorbing states and 5 stages. In this simple example, we have the following set of policy combinations: $\Upsilon = \{(0,0),(0,1),(1,0),(1,1)\}$, as we have 2 non-absorbing states. Each possible combination in $\Upsilon$ implies a policy that can be utilized in any stage. For example, $\upsilon=(1,0)$ means that we use the resource for state 0 but not for the ones in state 1,i.e., $\pi(0)=1,\pi(1)=0$. Furthermore, $\sigma^1,\sigma^2,$ and $\sigma^3$ all define a strategic path in $\mathcal{G}$, equivalently a strategy for the problem. Figure \ref{fig:exampleNetStrategy} also illustrates how $\sigma^1,\sigma^2,$ and $\sigma^3$ encode strategies.

Our ultimate goal for $\mathcal{G}$ is to find the strategy path $\sigma^*$ which provides the maximum total expected reward without violating the capacity constraints. To that end, we should specify the \textit{length} of an arc in a way that finding the feasible longest strategy path in $\mathcal{G}$ corresponds to finding $\Pi^*$. This is the point where we use Corollary \ref{cor:forward} which increases the efficiency by eliminating the iterations for $t'=1,\dotsc,t-1$ for the computation required for decision epoch $t$. Then, the length of arc $\iota=(\alpha', \alpha'') \in \Gamma$, denoted by $\eta_{\iota}$, is defined as follows:
\begin{itemize}
    \item if $\alpha' \in \phi_1$: 
    Then, $\eta_{\iota}$ is the sum of the expected rewards obtained in $t=1$ and $t=2$ based on the occupancy measures computed for $t=1$ (let's say $\mathcal{X}_1, \mathcal{Z}_1$) and $t=2$ (let's say $\mathcal{X}_2, \mathcal{Z}_2$). ($\mathcal{X}_1, \mathcal{Z}_1$) is computed based on $\mathbf{\theta}$ and $\pi^1$ implied by $\alpha'$ and ($\mathcal{X}_2, \mathcal{Z}_2$) is computed based on ($\mathcal{X}_1, \mathcal{Z}_1$) and $\pi^2$ implied by $\alpha''$. Then, we have the following expression:
    \begin{equation}
        \eta_{\iota} = \sum\limits_{\omega \in \Omega}\lambda_\omega\Big[(\mathcal{Z}^{\omega,1}+\mathcal{Z}^{\omega,2})R^D + \sum\limits_{i \in \tilde{\mathcal{S}}} (\mathcal{X}^{\omega1}_{i,\pi_i^1} r^\omega_{i,\pi_i^1} +\mathcal{X}^{\omega2}_{i,\pi_i^2} r^\omega_{i,\pi_i^2})\Big]
    \end{equation}
    
    \item if $\alpha' \in \phi_t: t \in \{2,\dotsc,T-3\}$: 
     Then, $\eta_{\iota}$ is the additional expected reward obtained by ($\mathcal{X}_{t+1}, \mathcal{Z}_{t+1}$) computed based on ($\mathcal{X}_{t}, \mathcal{Z}_{t}$) and $\pi^{t+1}$ implied by $\alpha''$. Then, we have the following expression:
    \begin{equation}
        \eta_{\iota} = \sum\limits_{\omega \in \Omega}\lambda_\omega \Big[\mathcal{Z}^{\omega,t+1}R^D + \sum\limits_{i \in \tilde{\mathcal{S}}} \mathcal{X}^{\omega,t+1}_{i, \pi_i^{t+1}}r^\omega_{i, \pi_i^{t+1}}\Big]
    \end{equation}
    
    \item if $\alpha' \in \phi_{T-2}$: 
    Then, $\eta_\iota$ is the sum of the additional expected reward obtained by ($\mathcal{X}_{T-1}, \mathcal{Z}_{T-1}$) computed based on ($\mathcal{X}_{T-2}, \mathcal{Z}_{T-2}$) and $\pi^{T-1}$ implied by $\alpha''$ and the additional expected reward obtained by ($\mathcal{Y}, \mathcal{Z}_{T}$) computed based on ($\mathcal{X}_{T-1}, \mathcal{Z}_{T-1}$). Then, we have the following expression:
    \begin{equation}
        \eta_{\iota} = \sum\limits_{\omega \in \Omega}\lambda_\omega\Big[(\mathcal{Z}^{\omega,T-1}+\mathcal{Z}^{\omega,T})R^D + \sum\limits_{i \in \tilde{\mathcal{S}}} \mathcal{X}^{\omega,T-1}_{i, \pi_i^{T-1}}r^\omega_{i, \pi_i^{T-1}}+\sum\limits_{i\in \tilde{\mathcal{S}}}Y^\omega_{i}R^\omega_i\Big]
    \end{equation}
\end{itemize}

It is important to notice that the length of an arc is not static but changes depending on the previous arc. With the appropriate length definition, our problem  is reduced to finding the feasible longest strategy path. In this regard, we develop a parallel approximate dynamic programming algorithm. In theory, the capacity constraints prevent us to claim the optimality. However, we will see that, through computational experiments, we generally get optimal solutions or very close solutions to global optima with the proposed algorithm illustrated in Algorithm \ref{Algorithm:ADP}.

\begin{algorithm}[ht!]
\caption{Approximate Dynamic Programming Algorithm (PADP)}
\begin{algorithmic}[1]
\scriptsize
\For{$\upsilon \in \Upsilon$}
        \State $\mathcal{X}^{\omega,1} \mathcal{Z}^{\omega,1} \gets$ \text{Policy Evaluation}($ \mathbf{\theta},\pi_t(\upsilon)$), $\forall \omega \in \Omega$
        \State $f_{\pi_1(\upsilon)}\gets 0$
        \If{$\exists \omega \in \Omega: n\sum\limits_{i \in \tilde{\mathcal{S}}}\mathcal{X}^{\omega,1}_{i,1}>C_1$}
        \State $\phi_1 \gets \phi_1 \setminus \{\pi_1(\upsilon)\} $
        \EndIf
\EndFor
\For{$t \in \{2,\dotsc,T-1\}$}
    \For{$\upsilon \in \Upsilon$}
        \State $f_{\pi_t(\upsilon)} \gets \underset{i \in \Xi_{\pi_{t}(\upsilon)}}{\max}\{f_i + \eta_{(i,\pi_{t}(\upsilon))}\}$
        \EndFor
 \EndFor
\Ensure $\{f_{\pi_t(\upsilon)}\}_{t \in \tilde{\mathcal{T}}, \upsilon \in \Upsilon}$
\end{algorithmic}\label{Algorithm:ADP}
\end{algorithm}

The algorithm tries to assign a \textit{value}, $f$, to each node in $\mathcal{G}$. $f_{\pi_t(\upsilon)}$ represents the value assigned to the node associated with stage $t$ and combination $\upsilon$. Basically, it shows the length of the feasible longest sub-strategic path found until the node $\pi_t(\upsilon) \in \phi_t$. The algorithm starts with the nodes in $\phi_1$. It computes $(\mathcal{X}_1, \mathcal{Z}_1)$ by using $\mathbf{\theta}$ and $\pi_1(\upsilon)$ for each $\upsilon \in \Upsilon$. The nodes which provide $\mathcal{X}_1$ occupancy measures in a way that the feasibility is violated are removed from $\phi_1$, and the value of the remaining nodes becomes zero. Then, iterations begin in a forward direction for the stages, i.e., $t=2,\dotsc,T-1$. In each stage $t$, all nodes within $\phi_t$ are visited. For each node $\pi_t(\upsilon) \in \phi_t$, the set of nodes that constitute feasible arcs targeting $\pi_t(\upsilon)$ form the set $\Xi_{\pi_t(\upsilon)}\subseteq \phi_{t-1}$. An arc is feasible if and only if the occupancy measures found for the target node don't violate the capacity constraints. Then, the value for $\pi_t(\upsilon)$, $f_{\pi_t(\upsilon)}$, is determined as $f_{\pi_t(\upsilon)} \gets \underset{i \in \Xi_{\pi_{t}(\upsilon)}}{\max}\{f_i + \eta_{(i,\pi_{t}(\upsilon))}\}$. The algorithm is terminated once all the nodes in $\phi_{T-1}$ are searched. Then, it returns the values found for each node in $\mathcal{G}$. The objective function value found by PADP algorithm is determined as the maximum value tackled by the nodes within $\phi_{T-1} $, i.e., $f^{\textit{PADP}} \triangleq \underset{i \in \phi_{T-1}}{\max}f_i$. Lastly, the algorithm performs very fast even if it is working with great number of scenarios and stages thanks to the parallel framework which tackles the scenarios separately.


\section{Case Study: Chronic Care Delivery Problem}
\label{section:case}

U.S. National Center for Health Statistics define a chronic condition as a disease which lasts three months or more, and they further claim that there is not any vaccine, medication, or cure that can stop it, or it does not disappear spontaneously. Therefore, it is assumed that recovery is not possible when a patient has a chronic condition. CDC demonstrate stroke, cancer, and diabetes as major chronic diseases. Bernell and Howard \cite{bernell2016use} further include hypertension, pulmonary conditions, and mental illnesses to the list of chronic diseases.

According to the records of National Health Insurance in the United States, 65\% of the population have multimorbidity, i.e., they have more than one chronic diseases. From societal perspective, a huge economic burden, approximately 80\% of Medicare spending, is caused by the patients with 4 or more chronic diseases \cite{wolff2002prevalence}. In addition, multimorbidity is more common in disadvantaged groups, which causes more social inequality \cite{britt2008prevalence}. Consequently, we are witnessing a changing paradigm called \textit{chronic care management} where polycronic patients' physical, mental and social needs are managed in one center. As in all systems, this approach has trade-offs since this system requires both fixed and variable costs. In this regard, the effective use of this resource is important. Therefore, the question of which patients should be targeted for chronic care delivery under capacity constraints becomes important. Here we will assume a hypothetical chronic care management system with a constrained resource. In this setting, the patient classes that will be eligible for the service for each period in the planning horizon should be determined at the system design stage.

The most distinctive and challenging part of designing a system for polychronic patients is that the underlying mechanism for the progress of health conditions may significantly vary from disease to disease. In the case of multimorbidity, the number of diseases, their severity, and which combinations of diseases occur can vary dramatically for each patient. For instance, a patient may have diabetes with type 2 and breast cancer, whereas another patient may have hypertension, asthma, and a cardiovascular disease. In this perspective, it is necessary to deal with transition probability and reward uncertainties, considering the range of chronic conditions, their breakdowns, and various levels of severity for the diseases. 

Patient targeting problem for chronic care delivery model is a suitable case study to test our model with its capacity constraints and inherently appeared parameter uncertainty within the system. In this regard, the associated MDP model is introduced in Section \ref{subsec:mdpCCD}.

\subsection{MDP Model}
\label{subsec:mdpCCD}

In this section, we conceptualize the dynamics of a patient in regard of chronic care delivery. In the rest of this section, we explain the components of the underlying MDP model.

\subsubsection{States}
\label{subsubsec:states}



A state involves two components. The first one is \textit{health status} which can be, without loss of generality, any measure indicating the physical well-being of a patient. In this study, we discretized this component into 3 categories: \textit{simple, moderate}, and \textit{complex}. The severity of conditions is getting worse from simple to complex.

The second component of a state includes the \textit{engagement} of a patient due to the statistics revealing that about 50$\%$ of premature deaths are caused by behaviors that could be changed \cite{national2015measuring}. Moreover, \cite{meng1999lifestyle} states that chronic diseases have very close link with lifestyle habits such as excessive alcohol use, poor nutrition, lack of physical activity, and tobacco use. In this regard, the term engagement corresponds to all behavioral aspects required to maintain patients’ own physical well-being. It is discretized by two alternatives: \textit{low} and \textit{high}. The engagement of patients increases from low to high.

In summary, we have six states coming from the pairwise combination of the physical and behavioral components: {Low-Simple (0), Low-Moderate (1), Low-Complex (2), High-Simple (3), High-Moderate (4), High-Complex (5)}. In addition to these states, we also include the state \textit{death (6)} as an absorbing state. 

\subsubsection{Actions}
\label{subsubsec:actions}

For each patient with a health status and a behavioral characteristic, an agent chooses either \textit{normal care} or \textit{special care}. Thus, the action space consists of two elements: 0, for normal care and 1, for special care. The action space is not binding for those reaching to the absorbing state.

\subsubsection{Transition Probabilities and Immediate Rewards}

Transition probabilities are subject to uncertainty due to highly combinatorial nature of multimorbidity. Hence, it is not possible to have data representing all possible transitions to obtain significant statistical inferences. Even if we have such data, the estimations are prone to large statistical errors because of limited or missing data for each possible transition. 
Nevertheless, it is possible to approximate the probability distributions by using expert opinions. We demonstrate one approach to employ these opinions. In this respect, we first generate hierarchical rules, i.e. relations, among transition probabilities based on the prior beliefs (see \ref{appendix:montecarlo_transitionprobs}). These are the rules that can be generated by any business analyst who is familiar with the system. Then, the base model the so-called \textit{nominal model} assumes that the parameters take the values the Monte Carlo algorithm depicted in Algorithm \ref{Algorithm:monte-carlo} converges to. The point where Monte Carlo converges may be seen as a center of the polyhedral set in which our nominal model for transition probabilities is likely to lie. Then, all scenarios are generated based on this center by adding uniform random noises around it. A similar approach is carried out for the immediate rewards with their own rules. Detailed description of scenario generation is explained in Section \ref{subsec:instanceGeneration}.


\begin{algorithm}[H]
\caption{Monte Carlo Approach for Transition Probability Estimation}
\begin{algorithmic}[1]
\scriptsize
\Require \text{number of iterations, expert opinions}

\State $\text{nIteration} \gets 1$ 
\State $\Bar{\mathcal{P}} \gets \emptyset$
\State $\Bar{\mathcal{Q}} \gets \emptyset$

\For{$\text{nIteration}\leq\text{number of iterations}$}
    \State $\mathcal{P}'' \gets $\text{A random realization which satisfies the required rules}
    \State $\mathcal{Q}'' \gets $\text{A random realization which satisfies the required rules}
    \State $\Bar{\mathcal{P}} \gets \Bar{\mathcal{P}} \cup \{\mathcal{P}''\}$
    \State $\Bar{\mathcal{Q}} \gets \Bar{\mathcal{Q}} \cup \{\mathcal{Q}''\}$
    \State $\text{nIteration} \gets \text{nIteration} + 1$ 
\EndFor
\State $\Hat{\mathcal{P}} \gets \text{mean}(\Bar{\mathcal{P}})$
\State $\Hat{\mathcal{Q}} \gets \text{mean}(\Bar{\mathcal{Q}})$
\Ensure $\Hat{\mathcal{P}}, \Hat{\mathcal{Q}}$
\end{algorithmic}\label{Algorithm:monte-carlo}
\end{algorithm}

Our Monte Carlo algorithm starts with empty set of patients. Then, a patient following the hierarchical rules is randomly generated in each iteration until reaching to the desired number of iterations. Then, the parameters of the nominal model are concluded as the mean of the corresponding parameters of the patients generated through the simulation.

\section{Computational Experiments}
\label{section:experiment}

In this section, we conduct computational experiments\footnote{All experiments are conducted on a computer with an Intel® Core ™ i7-6500U 2.39 GHz processor and 64 GB of RAM, with the Windows 10 operating system.} to obtain an understanding of (i) the contribution of the proposed model; (ii)  how effective our solution approach is. In this context, Section \ref{subsec:instanceGeneration} explains the way problem instances are generated. Then, Section \ref{subsec:impactVI} is devoted to investigate the impact of valid inequalities, whereas Section \ref{subsec:performanceADP} examines the computational aspects of the proposed algorithm. Next, Section \ref{subsec: valuePerfectInfo} analyzes the value of perfect information after Section \ref{subsec:valueStochastic} studies the value of stochastic solution. Finally, Section \ref{subsec:managerial_insights} provides managerial insights for potential practitioners with a special focus on the price of fairness, the value of flexibility, and capacity management.

\subsection{Instance Generation}
\label{subsec:instanceGeneration}

For chronic care delivery problem, it is not possible to have longitudinal data capturing all possible transitions in sufficient numbers. This is due to the large number of chronic conditions each of which has a variety of types. Furthermore, patients may have different numbers of chronic conditions of numerous severity. In this respect, it becomes an unrealistic expectation to have reliable time-series data of patients as in specialized services such as diabetes, HIV, and cancer treatments. Nevertheless, we have some prior beliefs about the ranks among the uncertain parameters. These are not complicated medical consequences, but simple justifiable facts allowing us to come up with a polyhedral uncertainty set for the nominal model. Point estimations, what is so-called nominal model, are generated through Monte Carlo approach depicted in Algorithm \ref{Algorithm:monte-carlo}. Then, different scenarios are generated by adding uniform noise to the nominal model. We refer to Hörman et al. (2013) and Zhang et al. (2017) for different scenario generation methodologies \cite{hormann2013automatic, zhang2017probabilistic}.

Let $\bar{\omega}_x$ denotes the nominal value of uncertain parameter $x$, whereas $\omega_x$ denotes the corresponding parameter value in scenario $\omega \in \Omega$. Given nominal model $\bar{\omega}$, scenario $\omega$ is generated as follows:
\begin{equation}
    \omega \triangleq \bigg\{\omega_x: \omega_x \sim \textit{UNIF}\Big((1-\epsilon)\bar{\omega}_x, (1+\epsilon)\bar{\omega}_x\Big)\bigg\}, \quad \epsilon \in (0,1)
\end{equation}
Then, elements of transition probability matrices are normalized so that they satisfy the conditions in (\ref{validstochastic})-(\ref{validstochastic2}).

A problem instance is uniquely characterized by 4 parameters. (i) $|\Omega|$, the number of scenarios in the solution sample; (ii) $T$, the number of decision epochs; (iii) $c=C/n$, the proportion of the population that can utilize the resource at a particular decision epoch; (iv) $\epsilon$, the maximum allowed variation from the nominal model. Without loss of generality, we assume that the capacities allocated per decision epoch are equal. Then, a problem instance with the given parameters is denoted by $\mathcal{I}(|\Omega|, T, c,\epsilon)$. Furthermore, prior probabilities are assumed to be equal. That is, the probability for a new patient to start the decision process in a state is equal for each non-absorbing state in $\tilde{\mathcal{S}}$, i.e., $\mathbb{P}\{s_1=i\}=\frac{1}{|\tilde{\mathcal{S}}|}, \forall i \in \tilde{\mathcal{S}}$. Lastly, the reward obtained after visiting the absorbing state, $r^D$, is determined as zero.

\subsection{Impact of Valid Inequalities}
\label{subsec:impactVI}

There are three components affecting the computational complexity of the problem: the number of decision epochs, scenarios, and states. In the case study, i.e., chronic care delivery problem, the size of the state space is fixed which is not problem since our objective is to embrace the problems with a large number of scenarios and/or stages. Therefore, we generated problem instances for a varying number of scenarios and decision epochs. To this end, we created 35 problem instances, each with the values of $ c = 0.4 $ and $ \epsilon = 0.25 $. While generating them, we use the set of different numbers of scenarios $\{5,10,25,50,100,250,500\}$ and the set of various numbers for decision epochs $\{5,10,20,30,40\}$.

In order to prepare an experiment analyzing the impact of the valid inequalities, we first solve (MIP-MMDP) by using the well-known commercial solver CPLEX 12.10 for all problem instances. Then, each valid inequality is separately included to the formulation, and solved again by the solver. All computational runs are imposed a time limit of 4 hours. In this time horizon, the solver could find the optimal solutions only for 21 instances among the 35 instances we generated. Table \ref{tab:VI} provides required solution times in seconds and relative percentage gap values for each instance and setting. We compare the required solution times for the 21 instances where the optimal solutions could be found, whereas the relative percentage gap values are compared for the rest of the instances. 

When we focus on the instances with zero gap values\footnote{Valid inequality 1 and 2 refer to Propositions \ref{VI:1} and \ref{VI:2}, respectively.}, we realize that valid inequality 1 provides the best solution time in $33.3 \%$ of the instances, whereas this statistic becomes $15.15\%$ for valid inequality 2. In the case where we look at the remaining instances, valid inequality 2 provides the best gap values for almost half of the instances, whereas valid inequality 1 provides it only for 2 problem instances. 

\begin{table}
    \centering
    \scriptsize
    \caption{Impact of Valid Inequalities in 35 Problem Instances}
    \label{tab:VI}
    \begin{tabular}{|rr|cc|cc|cc|}
\hline
           &            & \multicolumn{ 2}{|c}{(MIP-MMDP) Model} & \multicolumn{ 2}{|c}{Valid Inequality 1} & \multicolumn{ 2}{|c|}{Valid Inequality 2} \\
\hline
 nScenario &     nStage & Time (sec) &  Gap (\%) & Time (sec) &  Gap (\%) & Time (sec) &  Gap (\%) \\
\hline
\multicolumn{ 1}{|c}{5} &          5 &        0.2 &          0 &      0.247 &          0 &      0.965 &          0 \\

\multicolumn{ 1}{|c}{} &         10 &      2.371 &          0 &      \color{red}{2.36} &          0 &      2.538 &          0 \\

\multicolumn{ 1}{|c}{} &         20 &      8.241 &          0 &       8.95 &          0 &      \color{blue}{7.27} &          0 \\

\multicolumn{ 1}{|c}{} &         30 &      10.31 &          0 &       \color{red}{9.99} &          0 &     10.818 &          0 \\

\multicolumn{ 1}{|c}{} &         40 &       9.63 &          0 &      9.788 &          0 &     10.232 &          0 \\
\hline
\multicolumn{ 1}{|c}{10} &          5 &      0.225 &          0 &      0.221 &          0 &       \color{blue}{0.21} &          0 \\

\multicolumn{ 1}{|c}{} &         10 &        6.8 &          0 &      6.845 &          0 &      8.441 &          0 \\

\multicolumn{ 1}{|c}{} &         20 &     53.848 &          0 &     80.211 &          0 &    116.207 &          0 \\

\multicolumn{ 1}{|c}{} &         30 &    101.423 &          0 &     162.75 &          0 &    132.885 &          0 \\

\multicolumn{ 1}{|c}{} &         40 &    138.449 &          0 &    157.163 &          0 &    170.913 &          0 \\
\hline
\multicolumn{ 1}{|c}{25} &          5 &      0.511 &          0 &      0.532 &          0 &        \color{blue}{0.47} &          0 \\

\multicolumn{ 1}{|c}{} &         10 &    256.734 &          0 &      \color{red}{229.86} &          0 &    251.009 &          0 \\

\multicolumn{ 1}{|c}{} &         20 &   1891.068 &          0 &   1515.171 &          0 &    \color{blue}{1368.60} &          0 \\

\multicolumn{ 1}{|c}{} &         30 &   1425.692 &          0 &    \color{red}{1193.41} &          0 &   1783.103 &          0 \\

\multicolumn{ 1}{|c}{} &         40 &    3015.04 &          0 &   \color{red}{2818.66} &          0 &   3052.309 &          0 \\
\hline
\multicolumn{ 1}{|c}{50} &          5 &      9.664 &          0 &     \color{red}{9.486} &          0 &      9.846 &          0 \\

\multicolumn{ 1}{|c}{} &         10 &   3643.571 &          0 &   \color{red}{1176.563} &          0 &   2355.285 &          0 \\

\multicolumn{ 1}{|c}{} &         20 &  14401.456 &      0.6 &  14401.213 &      0.9 &  14401.537 &      0.8 \\

\multicolumn{ 1}{|c}{} &         30 &  14401.767 &      0.5 &  14439.982 &       0.01 &  14402.806 &      0.6 \\

\multicolumn{ 1}{|c}{} &         40 &   14470.37 &      1.1 &  14429.526 &      1.4 &  14430.794 &      1.2 \\
\hline
\multicolumn{ 1}{|c}{100} &          5 &    106.115 &          0 &    196.078 &          0 &    \color{blue}{96.072} &          0 \\

\multicolumn{ 1}{|c}{} &         10 &   9669.398 &          0 &  13778.265 &          0 &  14425.731 &      0.1 \\

\multicolumn{ 1}{|c}{} &         20 &  14404.549 &       2.0 &   14406.66 &      \color{red}{1.9} &  14409.823 &      2.3 \\

\multicolumn{ 1}{|c}{} &         30 &  14421.458 &      2.1 &  14406.278 &      2.1 &  14432.652 &      \color{blue}{1.9} \\

\multicolumn{ 1}{|c}{} &         40 &  14406.975 &      1.8 &  14434.273 &      2.1 &  14411.393 &      1.9 \\
\hline
\multicolumn{ 1}{|c}{250} &          5 &    193.304 &          0 &    339.143 &          0 &    272.681 &          0 \\

\multicolumn{ 1}{|c}{} &         10 &  14405.094 &      1.9 &  14407.678 &      1.8 &  14407.396 &      \color{blue}{1.7} \\

\multicolumn{ 1}{|c}{} &         20 &  14406.464 &      5.4 &  14407.621 &      5.4 &  14409.304 &      \color{blue}{4.8} \\

\multicolumn{ 1}{|c}{} &         30 &  14404.798 &      6.4 &  14405.341 &       6.0 &  14405.035 &      \color{blue}{5.3} \\

\multicolumn{ 1}{|c}{} &         40 &  14406.866 &      5.7 &  14407.235 &      6.1 &  14403.145 &      \color{blue}{5.5} \\
\hline
\multicolumn{ 1}{|c}{500} &          5 &    742.221 &          0 &    861.041 &          0 &   1103.164 &          0 \\

\multicolumn{ 1}{|c}{} &         10 &  14407.241 &      4.5 &  14402.248 &      5.4 &  14406.788 &      5.3 \\

\multicolumn{ 1}{|c}{} &         20 &  14406.308 &       7.0 &  14407.461 &      7.1 &  14410.242 &      \color{blue}{6.9} \\

\multicolumn{ 1}{|c}{} &         30 &  14406.978 &      7.7 &   14405.62 &      \color{red}{7.6} &   14406.41 &      7.7 \\

\multicolumn{ 1}{|c}{} &         40 &  14405.674 &      8.9 &  14409.774 &      8.9 &  14407.159 &      9.2 \\
\hline
\end{tabular}  
\end{table}   

\subsection{Computational Performance of the Parallel Approximate Dynamic Programming Algorithm}
\label{subsec:performanceADP}

Table \ref{tab:VI} shows that commercial solvers are not capable of solving problem instances with the desired number of scenarios for convergence. Therefore, an efficient algorithm providing solutions of high quality in a short span of time is essential for the applicability and the value of the developed model. In this context, a high quality solution refers to the ones that are optimal or very close to it and found within an acceptable amount of time.

In Section \ref{section:approach}, a parallel approximate dynamic programming algorithm is proposed as an alternative to solvers to provide an efficient algorithm required for the model. In this section, we test its capabilities in terms of time and optimality. In accordance with this purpose, we first use the 21 problem instances for which the solver could provide the optimal solutions in Section \ref{subsec:impactVI}. The proposed algorithm, PADP, is run for these problem instances, and the corresponding solution times in seconds and the percentage gap values are presented through Table \ref{tab:adp1}. PADP finds the optimal solutions for $42.857\%$ of the problem instances. For the remainder of them, the maximum deviation from the optimal solution is only $0.24\%$ which is quite acceptable for an algorithm trying to find an approximate solution. When we focus on the aggregate level performance metrics, we see an outstanding approximation so that the mean percentage gap value is only $0.073\%$. Another remarkable fact about the superiority of the algorithm is that it achieves these approximations in considerably short span of times. It finds the  solution about 1000 times faster than the solver.

In Table \ref{tab:adp2}, we compare the capabilities of CPLEX 12.10 and PADP for the instances the solver could not find the optimal solution in 4 hours in Section \ref{subsec:impactVI}. Recall that the solver loses its power for the instances with a large number of scenarios. The column \textit{improvement(\%)} shows the relative percentage improvement provided by PADP. A positive value in the column shows that the algorithm provides a better solution than the solver with 4 hours time limitation. On the contrary, the negative value indicates that the solution found by the solver in 4 hours provides a better objective function than the algorithm. Results suggest that the algorithm produces better results than the solver when the number of scenarios exceeds 250, i.e., for more realistic instances. The amount of improvement becomes more significant as the number of scenarios increases.

\begin{table}[p]
    \centering
    \scriptsize
    \caption{Computational Comparison of PADP and CPLEX 12.10 for Small-Medium Problem Instances}
    \begin{tabular}{|c|c|c|c|}
\hline
  Instance &  Time (MIP) & Time (PADP) &    GAP(\%) \\
\hline
         $\mathcal{I}(5,5,0.4,0.25)$ &        0.2 &      0.084 &          0 \\

         $\mathcal{I}(5,10,0.4,0.25)$ &      2.371 &      0.204 &          0 \\

         $\mathcal{I}(5,20,0.4,0.25)$ &      8.241 &      0.683 &       0.13 \\

         $\mathcal{I}(5,30,0.4,0.25)$ &      10.31 &      1.127 &       0.16 \\

         $\mathcal{I}(5,40,0.4,0.25)$ &       9.63 &      1.567 &       0.16 \\

         $\mathcal{I}(10,5,0.4,0.25)$ &      0.225 &      0.115 &          0 \\

         $\mathcal{I}(10,10,0.4,0.25)$ &        6.8 &       0.36 &          0 \\

         $\mathcal{I}(10,20,0.4,0.25)$ &     53.848 &      1.041 &       0.02 \\

         $\mathcal{I}(10,30,0.4,0.25)$ &    101.423 &      1.742 &       0.05 \\

         $\mathcal{I}(10,40,0.4,0.25)$ &    138.449 &      2.437 &       0.06 \\

         $\mathcal{I}(25,5,0.4,0.25)$ &      0.511 &      0.149 &          0 \\

         $\mathcal{I}(25,10,0.4,0.25)$ &    256.734 &      0.466 &          0 \\

         $\mathcal{I}(25,20,0.4,0.25)$ &   1891.068 &      1.332 &       0.14 \\

         $\mathcal{I}(25,30,0.4,0.25)$ &   1425.692 &      2.248 &        0.2 \\

         $\mathcal{I}(25,40,0.4,0.25)$ &    3015.04 &      3.176 &       0.22 \\

         $\mathcal{I}(50,5,0.4,0.25)$ &      9.664 &       0.17 &       0.02 \\

         $\mathcal{I}(50,10,0.4,0.25)$ &   3643.571 &      0.555 &       0.24 \\

         $\mathcal{I}(100,5,0.4,0.25)$ &    106.115 &      0.245 &          0 \\

         $\mathcal{I}(10,10,0.4,0.25)$ &   9669.398 &      0.802 &       0.13 \\

         $\mathcal{I}(250,5,0.4,0.25)$ &    193.304 &      0.304 &          0 \\

         $\mathcal{I}(500,5,0.4,0.25)$ &    742.221 &      0.434 &          0 \\
\hline
       min &        0.2 &      0.084 &          0 \\

       max &   9669.398 &      3.176 &       0.24 \\

      mean &   1013.563 &      0.916 &      0.073 \\
\hline
\end{tabular}  
    \label{tab:adp1}
\end{table}

\begin{table}[p]
    \centering
    \scriptsize
    \caption{Computational Comparison of PADP and CPLEX 12.10 for Large Problem Instances}
    \begin{tabular}{|c|c|c|c|c|}
\hline
  Instance & CPLEX Time(sec) & CPLEX Gap(\%) & PADP Time(sec) & Improvement(\%) \\
\hline
         $\mathcal{I}(50,20,0.4,0.25)$ &  14401.456 &        0.6 &      1.628 &      -0.33 \\

         $\mathcal{I}(50,30,0.4,0.25)$ &  14401.767 &        0.5 &      2.699 &      -0.36 \\

         $\mathcal{I}(50,40,0.4,0.25)$ &   14470.37 &        1.1 &      4.025 &      -0.37 \\

         $\mathcal{I}(100,20,0.4,0.25)$ &  14404.549 &          2 &      2.251 &       0.05 \\

         $\mathcal{I}(100,30,0.4,0.25)$ &  14421.458 &        2.1 &      3.635 &       0.01 \\

         $\mathcal{I}(100,40,0.4,0.25)$ &  14406.975 &        1.8 &        4.9 &      -0.14 \\

         $\mathcal{I}(250,10,0.4,0.25)$ &  14405.094 &        1.9 &       0.98 &       0.09 \\

         $\mathcal{I}(250,20,0.4,0.25)$ &  14406.464 &        5.4 &      3.005 &       0.66 \\

         $\mathcal{I}(250,30,0.4,0.25)$ &  14404.798 &        6.4 &      5.658 &       0.87 \\

         $\mathcal{I}(250,40,0.4,0.25)$ &  14406.866 &        5.7 &      7.724 &       0.13 \\

         $\mathcal{I}(500,10,0.4,0.25)$ &  14407.241 &        4.5 &      1.502 &       0.62 \\

         $\mathcal{I}(500,20,0.4,0.25)$ &  14406.308 &          7 &      5.662 &       1.13 \\

         $\mathcal{I}(500,30,0.4,0.25)$ &  14406.978 &        7.6 &      8.255 &       1.48 \\

         $\mathcal{I}(500,40,0.4,0.25)$ &  14405.674 &        8.9 &     10.203 &       2.66 \\
\hline
\end{tabular}  

\label{tab:adp2}
\end{table}

\subsection{Value of Stochastic Solution}
\label{subsec:valueStochastic}

If agents do not take parameter uncertainty into account then the actions will become sub-optimal. With this in mind, in this section, we investigate the additional value of incorporating stochastic solution. In other words, we estimate the potential loss in objective function value, if we don't consider the transition probability and reward uncertainties. There is certainly a strong relation between the value of stochastic solution and the way we generate the problem instances. The main assumption we make during the generation of data is the maximum variation from the nominal model. Thus, we conduct sensitivity analysis for different values of $\epsilon$. Particularly, we choose three values of $\epsilon$ for the experiments: $0.10, 0.25,$ and $0.50$. 

We also postulate uniform distribution when we consider the deviations from the nominal model, and by force of uniform distribution, it is not possible to observe any realization for the parameters outside of the sphere determined by $\epsilon$. However, observing lots of outliers which lie in the exterior region of the sphere is very likely, considering the large number of combinations within the nature of the problem. Thus, the problem includes more randomness than our problem instances. Neverthless, there are two advantages of this approach. Firstly, the value of stochastic solution asserted for these problem instances intuitively provides a lower bound for the actual value since we don't allow outliers. Therefore, if we find a considerable value of stochastic solution then we can assert that the actual value is expected to perform higher than it, which makes our contribution realistic. The second advantage is that this way of generating data materializes itself without loss of generality. That is, results associated with other experiments about the complexity of the model and the efficiency of our algorithm are not effected by the assumptions made during the generation of the problem instances.

Expected value of stochastic solution, EVSS($\%$), under solution sample $\Omega=\{1,\dotsc,\omega,\dotsc,|\Omega|\}$ is computed as follows:

\begin{equation}
    \label{eq:value}
    \text{EVSS}(\%)\triangleq \frac{1}{|\Omega|}\sum_{\omega \in \Omega}\frac{f^{\textit{MIP-MMDP}}-f^\omega(\Omega)}{f^\omega(\Omega)}\times 100
\end{equation}
where $f^{\textit{MIP-MMDP}}$ is the total expected reward of the optimal strategy found for $\Omega$, i.e., the optimal objective function value of (MIP-MMDP) for $\Omega$, and $f^\omega(\Omega)$ is the total expected reward of following the optimal strategy for scenario $\omega$, where the realization of nature is still explained by all scenarios in $\Omega$. It is clear that $f^\omega(\Omega)$ yields a lower bound for $f^{\textit{MIP-MMDP}}$, i.e., $f^\omega(\Omega) \leq f^{\textit{MIP-MMDP}}, \forall \omega \in \Omega$ due to the sub-optimal behavior of the strategy governed by a single scenario. Thus, $\frac{f^{\textit{MIP-MMDP}}-f^\omega(\Omega)}{f^\omega(\Omega)}$ shows the relative loss in the objective function value if agents consider only scenario $\omega \in \Omega$ as a deterministic model. Then, $EVSS(\%)$ gives the expected value of stochastic solution. \ref{appendix:evaluate_f} shows how $f^\omega(\Omega)$ is evaluated for scenario $\omega \in \Omega$.

\begin{table}[H]
    \centering
    \scriptsize
    \caption{Expected Value of Stochastic Solution for Different $\epsilon$ and $T$ Values with 200 Scenarios}
\begin{tabular}{|c|c|c|}
\hline
   $\epsilon$ &     nStage & Value of stochastic solution($\%$) \\
\hline
\multicolumn{ 1}{|c|}{0.1} &       5 &       2.97 \\

\multicolumn{ 1}{|c|}{} &         10 &      4.52 \\

\multicolumn{ 1}{|c|}{} &         20 &      4.47 \\

\multicolumn{ 1}{|c|}{} &         30 &      4.53 \\
\hline
\multicolumn{ 1}{|c|}{0.25} &      5 &      4.03 \\

\multicolumn{ 1}{|c|}{} &         10 &      5.70 \\

\multicolumn{ 1}{|c|}{} &         20 &      6.83 \\

\multicolumn{ 1}{|c|}{} &         30 &      7.31 \\
\hline
\multicolumn{ 1}{|c|}{0.5} &       5 &      3.95 \\

\multicolumn{ 1}{|c|}{} &         10 &      5.89 \\

\multicolumn{ 1}{|c|}{} &         20 &      7.80 \\

\multicolumn{ 1}{|c|}{} &         30 &      9.30 \\
\hline
\end{tabular}  
\label{tab:value}
\end{table}

Table \ref{tab:value} shows the expected value of stochastic solution for different values of $\epsilon$ and $T$ when the number of scenarios is 200. Results suggest that there is a significant value of incorporating stochastic solution. Expected value is increasing as the value of $\epsilon$ increases. Another important observation is that it also increases with the growing number of decision epochs which implies that our approach matters more for long-term problems.

\subsection{Value of Perfect Information}
\label{subsec: valuePerfectInfo}

Value of perfect information (VPI) is defined as the maximum cost a policy-maker is willing to pay in order to obtain the perfect information about the uncertain parameters. Since we are dealing with a maximization problem, the objective function of single scenario deterministic problem, denoted by $f^\omega$, becomes an upper bound for $f^{\text{MIP-MMDP}}$, i.e., $f^{\text{MIP-MMDP}} \leq f^\omega, \forall \omega \in \Omega$. Then, the expected VPI (EVPI) under the set of scenarios $\Omega$ is computed as follows:
\begin{equation}
    \text{EVPI} = \frac{\sum\limits_{\omega \in \Omega}f^\omega}{|\Omega|} - f^{\text{MIP-MMDP}}
\end{equation}

Then, the percentage expected VPI, denoted by EVPI(\%) is defined as $\frac{\text{EVPI}}{f^{\text{MIP-MMDP}}} \times 100$. Table \ref{tab:EVPI} illustrates EVPI(\%) for different values of $\epsilon$ and $T$ when the solution sample consists of 200 scenarios. We compute EVPI(\%) for each pairwise combinations of $\epsilon=0.1, 0.25, 0.5$ and $T=5,10,20,30$. Results suggest that the average EVPI(\%) value is $7.88 \%$. They also indicate that EVPI(\%) increases with the increasing number of $\epsilon$ and $T$.

\begin{table}[ht]
    \centering
    \scriptsize
    \caption{Expected Value of Perfect Information for Different $\epsilon$ and $T$ Values with 200 Scenarios}
\begin{tabular}{|c|c|c|}
\hline
   $\epsilon$ &     nStage & Value of perfect information($\%$) \\
\hline
\multicolumn{ 1}{|c|}{0.1} &       5 &      1.72 \\

\multicolumn{ 1}{|c|}{} &         10 &      2.14 \\

\multicolumn{ 1}{|c|}{} &         20 &      1.84 \\

\multicolumn{ 1}{|c|}{} &         30 &      1.76 \\
\hline
\multicolumn{ 1}{|c|}{0.25} &      5 &      5.62 \\

\multicolumn{ 1}{|c|}{} &         10 &      7.09 \\

\multicolumn{ 1}{|c|}{} &         20 &      5.81 \\

\multicolumn{ 1}{|c|}{} &         30 &      5.17 \\
\hline
\multicolumn{ 1}{|c|}{0.5} &       5 &      12.44 \\

\multicolumn{ 1}{|c|}{} &         10 &      16.36 \\

\multicolumn{ 1}{|c|}{} &         20 &      17.66 \\

\multicolumn{ 1}{|c|}{} &         30 &      16.92 \\
\hline
\end{tabular}  
\label{tab:EVPI}
\end{table}

\subsection{Managerial Insights}
\label{subsec:managerial_insights}

In this section, we investigate managerial aspects of applying (MIP-MMDP). In this context, we first examine the price of fairness, that is, the loss in total expected reward caused by ethical issues. Then, we focus on the value of having flexibility in operations. Lastly, we analyze the impact of different capacity levels . It is important to notice that each part, except the one measuring the effect of different capacity levels, requires certain modifications in the model. For these parts, we change the required components of the model, and then solve them with the solver. That is, we do not modify our algorithm for the new cases, which is beyond the scope of our study. Thus, the analyzes are limited with the problem instances solved by the solver in 4 hours, i.e., small-medium problem instances. However, we think that even this can be useful to understand the potential values and drawbacks of the model for practical issues.

\subsubsection{Price of Fairness}
\label{subsubsec:price_of_fairness}
(MIP-MMDP) is limited to deterministic policies, even a randomized policy can promise a better objective function value. Nevertheless, we require deterministic policies due to the fairness concerns raised in healthcare operations \cite{cevik2018analysis, merakli2019risk}. This is also the case for other domains such as humanitarian relief operations \cite{steimle2018multi}. The term fairness relates to the ethical issues which oblige policy-makers to follow the same policy for all patients in the same state. Hence, (MIP-MMDP) with randomized policies, i.e., without (\ref{RNEF_deterministicPol1}) - (\ref{RNEF_domain4}), becomes a relaxation of the original formulation. Then, \textit{price of fairness (\%)}, determined as $\frac{f^r-f^{\text{MIP-MMDP}}}{f^r} \times 100$, corresponds to the loss in the objective function value because of the constraints which eliminate randomized policies, where $f^r$ denotes the optimal objective function value of the model allowing randomized policies.

\begin{table}[h]
    \centering
    \scriptsize
    \caption{Price of Fairness for Small-Medium Problem Instances}
\begin{tabular}{|ccc|ccc|}
\hline
{nScenario} & {nStage} & Price of Fairness (\%) & {nScenario} & {nStage} & Price of Fairness (\%) \\
\hline
5 &  5 & 4.10 & 25 & 5 & 5.96\\

5 & 10 & 4.66 & 25 & 10 &  6.98\\

5 & 20 & 4.15 & 25 & 20 &  6.45\\

5 & 30 & 3.89 & 25 & 30 &  6.31\\

5 & 40 & 3.80 & 25 & 40 &  6.27\\

10 & 5 & 4.19 & 50 & 5 &  6.81\\

10 & 10 &  5.02 & 50 & 10 &  7.99\\

10 & 20 &  4.89 & 100 & 5 & 7.13\\

10 & 30 &  4.82 & 100 & 10 & 8.39\\

10 & 40 &  4.80 & 250 & 5 & 7.08\\
& & & 500 & 5 & 7.07 \\
\hline
\end{tabular}
\label{tab:price_fairness}
\end{table}

Price of fairness (\%) for each small-medium problem instance is demonstrated in Table \ref{tab:price_fairness}. Results point out that the average loss caused by the ethical concerns is $5.75\%$. It stems from the fact that a portion of the capacity at each decision epoch becomes idle due to the deterministic policy constraints. Allowing randomized policies enables the model to fulfill the capacity, even if only a subset of patients in a state utilize the resource. Results also show that as the number of scenarios increases, price of fairness tends to increase. 

\subsubsection{Value of Flexibility}
\label{subsec:operational_flexibility}

We have so far assumed flexibility in actions, where in each period the actions are independent of the actions taken in other periods. That is, one state may be accepted in stage $t$ and not accepted in stages $t' \in \tilde{\mathcal{T}}\setminus\{t\}$. In practice, policy makers may not have this flexibility; a state is either entitled to the resource for all decision epochs from the beginning or not at all. Therefore, here we consider the case in which once action is 1 for a state in a particular decision epoch, then it must also be 1 for other decision epochs. Mathematically, this corresponds to adding the following constraint to (MIP-MMDP): $\pi_i^t - \pi_i^{t-1}=0, \forall i \in \tilde{\mathcal{S}}, t \in \tilde{\mathcal{T}}\setminus\{1\}$. It is clear that including it narrows the feasible region so that the previous objective function value becomes an upper bound for the new model, i.e., $\tilde{f} \leq f^{\text{MIP-MMDP}}$ where $\tilde{f}$ is the optimal objective function value under the same policy approach. In this study, flexibility implies that policy-makers can change their course of actions over time. Since (MIP-MMDP) allows it, \textit{value of flexibility (\%)} is determined as $\frac{f^{\text{MIP-MMDP}}-\tilde{f}}{\tilde{f}} \times 100$.

\begin{table}
    \centering
    \scriptsize
    \caption{Value of Flexibility for Small-Medium Problem Instances}
\begin{tabular}{|ccc|ccc|}
\hline
{nScenario} & {nStage} & value of Flexibility (\%) & {nScenario} & {nStage} & value of Flexibility (\%) \\
\hline
5 &  5 & 2.53 & 25 & 5 & 2.12\\

5 & 10 & 5.38 & 25 & 10 &  5.34\\

5 & 20 & 10.24 & 25 & 20 &  13.33\\

5 & 30 & 12.39 & 25 & 30 &  17.61\\

5 & 40 & 13.25 & 25 & 40 &  19.59\\

10 & 5 & 2.16 & 50 & 5 &  1.55\\

10 & 10 &  5.88 & 50 & 10 &  4.73\\

10 & 20 &  11.81 & 100 & 5 & 1.23\\

10 & 30 &  15.11 & 100 & 10 & 4.15\\

10 & 40 &  16.63 & 250 & 5 & 1.21\\
& & & 500 & 5 & 1.25 \\
\hline
\end{tabular}
\label{tab:value_flexibility}
\end{table}

\begin{figure}[ht]
\caption{Average Value of Flexibility for $T=5,10,20,30,40$}
  \includegraphics[width=\linewidth]{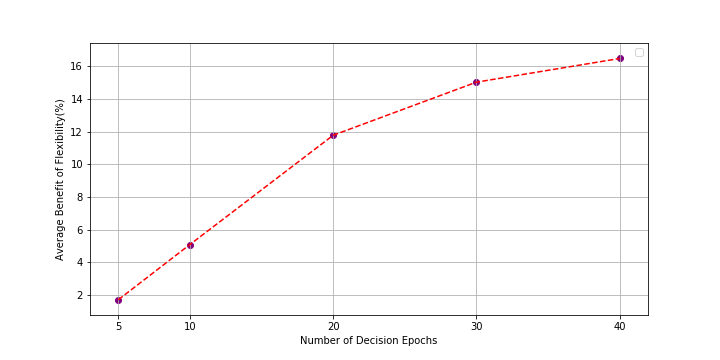}
  \label{fig:value_flexibility}
\end{figure}

Value of flexibility (\%) is computed for each small-medium problem instance, and the values are represented in Table \ref{tab:value_flexibility}. It shows that $7.98\%$ additional value can be obtained by (MIP-MMDP) instead of using the same policy approach. In other words, embracing flexibility yields almost $7.98 \%$ additional total expected reward. Furthermore, average values for $T=5,10,20,30,$ and $40$ are computed and summarized in Figure \ref{fig:value_flexibility}. It demonstrates that as number of decision epochs increases, the value of embracing flexibility also increases dramatically since the suboptimal behavior is retained for longer periods. This also shows the relative importance of our model for long-term planning problems compared to the ones for short-terms.

\subsubsection{Capacity Management}
\label{subsubsec:capacity_management}

\begin{figure}[H]
\caption{Total Expected Rewards under Different Stage Capacities and Decision Epochs for the Model Restricted with Deterministic Policies}
\includegraphics[width=\linewidth]{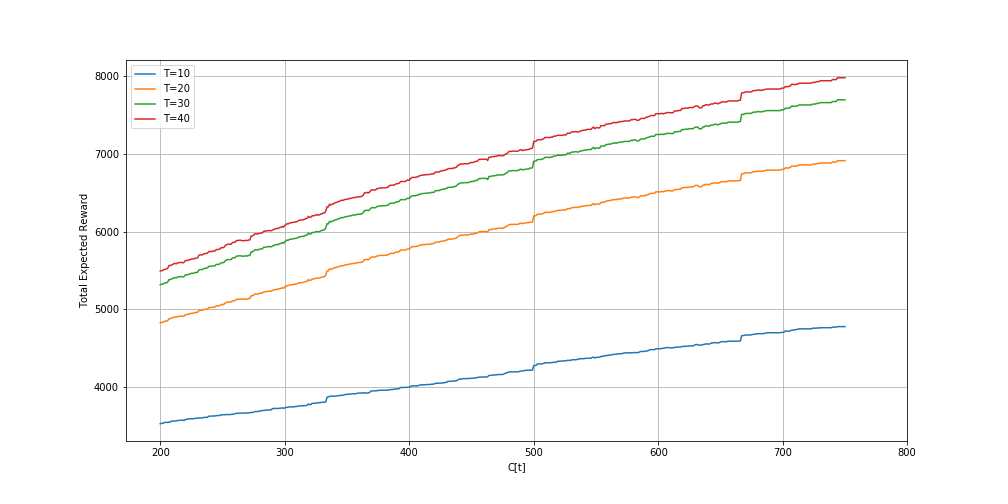}
\label{fig:capacity_management}
\end{figure}

\begin{figure}[H]
\caption{Total Expected Rewards under Different Stage Capacities and Decision Epochs for the Model  with Randomized Policies}
\includegraphics[width=\linewidth]{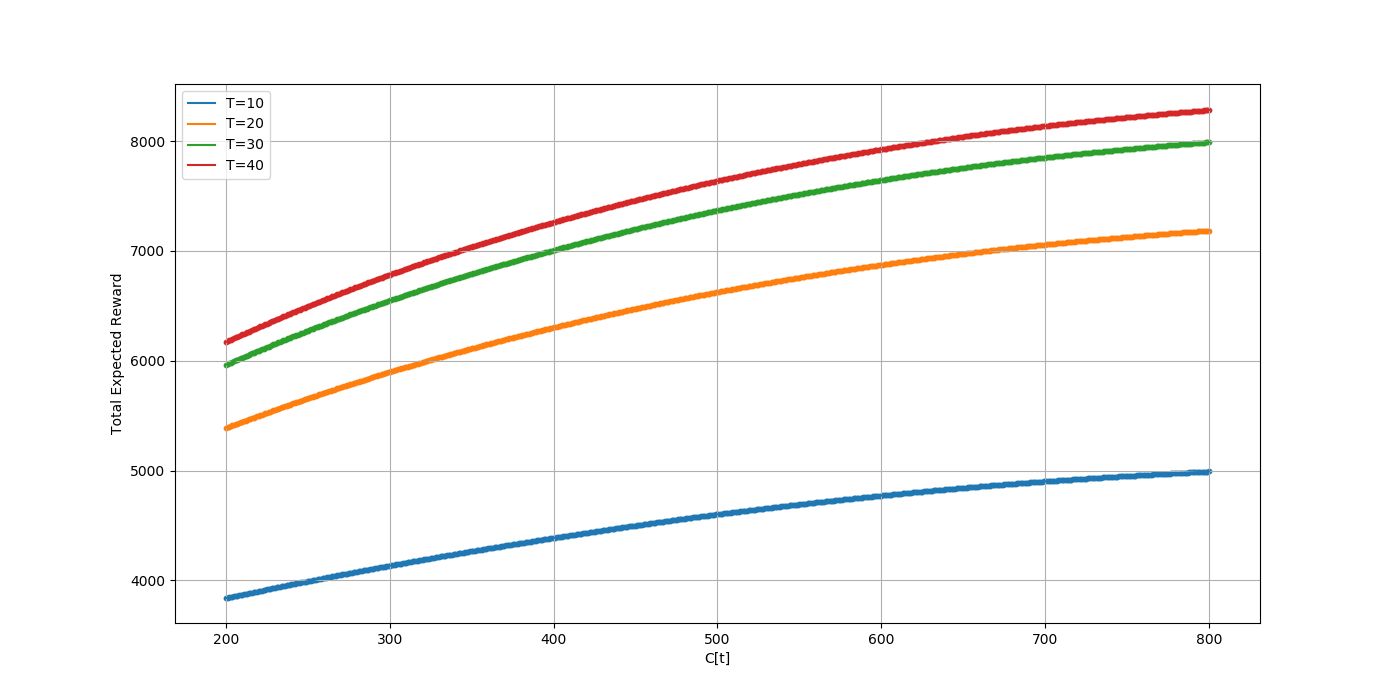}
\label{fig:capacity_management_randomized}
\end{figure}

In this section, we analyze the additional total expected reward that can be obtained by adding one more capacity at each period. To this end, we find total expected rewards for each instance from $\bigcup\limits_{c=0.2}^{0.8} \mathcal{I}(200,t,c,0.25)$ for $t \in \{10,20,30,40\}$. In words, for $t \in \{10,20,30,40\}$, we start with stage capacity 200 and gradually increases it until 800 for a fixed solution sample which consists of 200 scenarios. We also follow this procedure for the case where randomized policies are allowed. Figure \ref{fig:capacity_management} and Figure \ref{fig:capacity_management_randomized} illustrate the results of the experiment for the cases with deterministic policies and randomized policies, respectively. The main difference between deterministic and randomized policies become more distinguishable when we compare the graphs in Figure \ref{fig:capacity_management} and Figure \ref{fig:capacity_management_randomized}. 
Table \ref{fig:capacity_management_randomized} clearly shows that the graph tends to perform concave characteristics. That is, the marginal utility of adding one unit of stage capacity is decreasing, as the value of stage capacity increases. Unlike the graph in Figure $\ref{fig:capacity_management_randomized}$, we do not observe smooth increases within the graph in Figure \ref{fig:capacity_management}. That is, increasing the capacity for deterministic policies does not necessarily improves the current conditions since the model still may not find an available group to fulfill the remaining capacity. The frequently observed constant lines in Figure \ref{fig:capacity_management} corresponds to these cases. Thus, system managers embracing deterministic policies must take this situation into account when they need to decide how much to increase the capacity.

\section{Conclusion}
\label{section:conclusion}

In this study, we introduce a multi-model Markov decision process model with capacity constraints, motivated from medical resource allocation problems. Multi-model approach implies relaxation of the traditional assumption which ensures deterministic transition probabilities and rewards. Capacity constraints pose limitations on the selection of particular actions. To the best of our knowledge, the proposed model is the first constrained MDP model considering the parameter uncertainty. 

We develop a mixed integer programming (MIP) formulation for the model in order to find the strategy which maximizes the total expected reward without violating capacity constraints. It corresponds to an extensive-form formulation for the underlying
two-stage stochastic integer program. The scenarios used in the extensive-form formulation basically corresponds to different MDP models. The MIP formulation becomes easily intractable as the size of the problem grows. In this regard, we propose a parallel approximate dynamic programming algorithm leveraging the problem structure. We also propose two valid inequalities in the hope that they can strengthen the formulation.

We test test our model and algorithm using a chronic care delivery problem as an example. We adopt a Monte Carlo approach which uses prior beliefs about system dynamics to generate problem instances. Extensive computational experiments are then performed to test the computational aspects of the proposed algorithm as well as the value of our model. We show that the algorithm works pretty fast so that it generates solutions in seconds even for very large problem instances. It also offers very high quality solutions in terms of optimality. However, while our algorithm embraces the cases with a large number of scenarios and stages, it does not address the last level of complexity which is the state space. 

There are several potential approaches to complement our work. First, a new algorithm can be proposed that also works fast for large state space problems. In this context, it may be promising to address constrained reinforcement learning models. Secondly, for some medical decision making problems such as breast cancer screening, we may not have have perfect information of the health status of patient \cite{cevik2018analysis}. To deal with these kind of challenges, our model can be extended for POMDP models.

\newpage
\setstretch{1.0}
\bibliographystyle{apalike}
\bibliography{references}

\begin{thebibliography}{}

\bibitem[Alagoz et~al., 2010]{alagoz2010markov}
Alagoz, O., Hsu, H., Schaefer, A.~J., and Roberts, M.~S. (2010).
\newblock Markov decision processes: a tool for sequential decision making
  under uncertainty.
\newblock {\em Medical Decision Making}, 30(4):474--483.

\bibitem[Alagoz et~al., 2004]{alagoz2004optimal}
Alagoz, O., Maillart, L.~M., Schaefer, A.~J., and Roberts, M.~S. (2004).
\newblock The optimal timing of living-donor liver transplantation.
\newblock {\em Management Science}, 50(10):1420--1430.

\bibitem[Altman, 1999]{altman1999constrained}
Altman, E. (1999).
\newblock {\em Constrained Markov decision processes}, volume~7.
\newblock CRC Press.

\bibitem[Ayvaci et~al., 2012]{ayvaci2012effect}
Ayvaci, M.~U., Alagoz, O., and Burnside, E.~S. (2012).
\newblock The effect of budgetary restrictions on breast cancer diagnostic
  decisions.
\newblock {\em Manufacturing \& Service Operations Management}, 14(4):600--617.

\bibitem[Ben-Tal et~al., 2009]{ben2009robust}
Ben-Tal, A., El~Ghaoui, L., and Nemirovski, A. (2009).
\newblock {\em Robust optimization}, volume~28.
\newblock Princeton University Press.

\bibitem[Ben-Tal and Nemirovski, 2002]{ben2002robust}
Ben-Tal, A. and Nemirovski, A. (2002).
\newblock Robust optimization--methodology and applications.
\newblock {\em Mathematical programming}, 92(3):453--480.

\bibitem[Bernell and Howard, 2016]{bernell2016use}
Bernell, S. and Howard, S.~W. (2016).
\newblock Use your words carefully: what is a chronic disease?
\newblock {\em Frontiers in public health}, 4:159.

\bibitem[Bertsimas et~al., 2011]{bertsimas2011theory}
Bertsimas, D., Brown, D.~B., and Caramanis, C. (2011).
\newblock Theory and applications of robust optimization.
\newblock {\em SIAM review}, 53(3):464--501.

\bibitem[Boucherie and Van~Dijk, 2017]{boucherie2017markov}
Boucherie, R.~J. and Van~Dijk, N.~M. (2017).
\newblock {\em Markov decision processes in practice}, volume 248.
\newblock Springer.

\bibitem[Britt et~al., 2008]{britt2008prevalence}
Britt, H.~C., Harrison, C.~M., Miller, G.~C., and Knox, S.~A. (2008).
\newblock Prevalence and patterns of multimorbidity in australia.
\newblock {\em Medical Journal of Australia}, 189(2):72--77.

\bibitem[Cevik et~al., 2018]{cevik2018analysis}
Cevik, M., Ayer, T., Alagoz, O., and Sprague, B.~L. (2018).
\newblock Analysis of mammography screening policies under resource
  constraints.
\newblock {\em Production and Operations Management}, 27(5):949--972.

\bibitem[Council et~al., 2015]{national2015measuring}
Council, N.~R., on~Population, C., et~al. (2015).
\newblock {\em Measuring the Risks and Causes of Premature Death: Summary of
  Workshops}.
\newblock National Academies Press.

\bibitem[Denton, 2018]{denton2018optimization}
Denton, B.~T. (2018).
\newblock Optimization of sequential decision making for chronic diseases: From
  data to decisions.
\newblock In {\em Recent Advances in Optimization and Modeling of Contemporary
  Problems}, pages 316--348. INFORMS.

\bibitem[Denton et~al., 2009]{denton2009optimizing}
Denton, B.~T., Kurt, M., Shah, N.~D., Bryant, S.~C., and Smith, S.~A. (2009).
\newblock Optimizing the start time of statin therapy for patients with
  diabetes.
\newblock {\em Medical Decision Making}, 29(3):351--367.

\bibitem[Deo et~al., 2013]{deo2013improving}
Deo, S., Iravani, S., Jiang, T., Smilowitz, K., and Samuelson, S. (2013).
\newblock Improving health outcomes through better capacity allocation in a
  community-based chronic care model.
\newblock {\em Operations Research}, 61(6):1277--1294.

\bibitem[Givan et~al., 2000]{givan2000bounded}
Givan, R., Leach, S., and Dean, T. (2000).
\newblock Bounded-parameter markov decision processes.
\newblock {\em Artificial Intelligence}, 122(1-2):71--109.

\bibitem[Gorissen et~al., 2015]{gorissen2015practical}
Gorissen, B.~L., Yan{\i}ko{\u{g}}lu, {\.I}., and den Hertog, D. (2015).
\newblock A practical guide to robust optimization.
\newblock {\em Omega}, 53:124--137.

\bibitem[Grand-Clement et~al., 2020]{grand2020robust}
Grand-Clement, J., Chan, C.~W., Goyal, V., and Escobar, G. (2020).
\newblock Robust policies for proactive icu transfers.
\newblock {\em arXiv preprint arXiv:2002.06247}.

\bibitem[H{\"o}rmann et~al., 2013]{hormann2013automatic}
H{\"o}rmann, W., Leydold, J., and Derflinger, G. (2013).
\newblock {\em Automatic nonuniform random variate generation}.
\newblock Springer Science \& Business Media.

\bibitem[Iyengar, 2005]{iyengar2005robust}
Iyengar, G.~N. (2005).
\newblock Robust dynamic programming.
\newblock {\em Mathematics of Operations Research}, 30(2):257--280.

\bibitem[Mannor et~al., 2016]{mannor2016robust}
Mannor, S., Mebel, O., and Xu, H. (2016).
\newblock Robust mdps with k-rectangular uncertainty.
\newblock {\em Mathematics of Operations Research}, 41(4):1484--1509.

\bibitem[Mannor et~al., 2007]{mannor2007bias}
Mannor, S., Simester, D., Sun, P., and Tsitsiklis, J.~N. (2007).
\newblock Bias and variance approximation in value function estimates.
\newblock {\em Management Science}, 53(2):308--322.

\bibitem[Meng et~al., 1999]{meng1999lifestyle}
Meng, L., Maskarinec, G., Lee, J., and Kolonel, L.~N. (1999).
\newblock Lifestyle factors and chronic diseases: application of a composite
  risk index.
\newblock {\em Preventive medicine}, 29(4):296--304.

\bibitem[Merakl{\i} and K{\"u}{\c{c}}{\"u}kyavuz, 2019]{merakli2019risk}
Merakl{\i}, M. and K{\"u}{\c{c}}{\"u}kyavuz, S. (2019).
\newblock Risk aversion to parameter uncertainty in markov decision processes
  with an application to slow-onset disaster relief.
\newblock {\em IISE Transactions}, pages 1--21.

\bibitem[Nilim and El~Ghaoui, 2004]{nilim2004robustness}
Nilim, A. and El~Ghaoui, L. (2004).
\newblock Robustness in markov decision problems with uncertain transition
  matrices.
\newblock In {\em Advances in neural information processing systems}, pages
  839--846.

\bibitem[Puterman, 2014]{puterman2014markov}
Puterman, M.~L. (2014).
\newblock {\em Markov decision processes: discrete stochastic dynamic
  programming}.
\newblock John Wiley \& Sons.

\bibitem[Satia and Lave~Jr, 1973]{satia1973markovian}
Satia, J.~K. and Lave~Jr, R.~E. (1973).
\newblock Markovian decision processes with uncertain transition probabilities.
\newblock {\em Operations Research}, 21(3):728--740.

\bibitem[Seifert et~al., 2016]{seifert2016dynamic}
Seifert, R.~W., Tancrez, J.-S., and Bi{\c{c}}er, I. (2016).
\newblock Dynamic product portfolio management with life cycle considerations.
\newblock {\em International Journal of Production Economics}, 171:71--83.

\bibitem[Shechter et~al., 2008]{shechter2008optimal}
Shechter, S.~M., Bailey, M.~D., Schaefer, A.~J., and Roberts, M.~S. (2008).
\newblock The optimal time to initiate hiv therapy under ordered health states.
\newblock {\em Operations Research}, 56(1):20--33.

\bibitem[Sinha and Ghate, 2016]{sinha2016policy}
Sinha, S. and Ghate, A. (2016).
\newblock Policy iteration for robust nonstationary markov decision processes.
\newblock {\em Optimization Letters}, 10(8):1613--1628.

\bibitem[Steimle et~al., 2018]{steimle2018multi}
Steimle, L.~N., Kaufman, D.~L., and Denton, B.~T. (2018).
\newblock Multi-model markov decision processes.
\newblock {\em Optimization Online URL http://www. optimization-online.
  org/DB\_FILE/2018/01/6434. pdf}.

\bibitem[Tewari and Bartlett, 2007]{tewari2007bounded}
Tewari, A. and Bartlett, P.~L. (2007).
\newblock Bounded parameter markov decision processes with average reward
  criterion.
\newblock In {\em International Conference on Computational Learning Theory},
  pages 263--277. Springer.

\bibitem[White~III and Eldeib, 1994]{white1994markov}
White~III, C.~C. and Eldeib, H.~K. (1994).
\newblock Markov decision processes with imprecise transition probabilities.
\newblock {\em Operations Research}, 42(4):739--749.

\bibitem[Wiesemann et~al., 2013]{wiesemann2013robust}
Wiesemann, W., Kuhn, D., and Rustem, B. (2013).
\newblock Robust markov decision processes.
\newblock {\em Mathematics of Operations Research}, 38(1):153--183.

\bibitem[Wolff et~al., 2002]{wolff2002prevalence}
Wolff, J.~L., Starfield, B., and Anderson, G. (2002).
\newblock Prevalence, expenditures, and complications of multiple chronic
  conditions in the elderly.
\newblock {\em Archives of internal medicine}, 162(20):2269--2276.

\bibitem[Zhang et~al., 2017a]{zhang2017robust}
Zhang, Y., Steimle, L., and Denton, B. (2017a).
\newblock Robust markov decision processes for medical treatment decisions.
\newblock {\em Optimization online}.

\bibitem[Zhang et~al., 2017b]{zhang2017probabilistic}
Zhang, Y., Wu, H., Denton, B.~T., Wilson, J.~R., and Lobo, J.~M. (2017b).
\newblock Probabilistic sensitivity analysis on markov decision processes with
  uncertain transition probabilities: an application in evaluating treatment
  decisions for type 2 diabetes.

\end{thebibliography}

\newpage
\appendix

\section{Proof of Valid Inequalities}
\label{appendix:proof_valid}

\subsection{Proof of Proposition \ref{VI:1}}

\begin{itemize}
\item for $t=1$:

    We are given that $\sum\limits_{a \in \mathcal{A}}\mathcal{X}^{\omega,1}_{ia}=\theta_i$ for each $i \in \tilde{\mathcal{S}}$ and $\omega \in \Omega$ by (\ref{RNEF_flow1}). Then,
    \begin{align*}
        \sum_{i \in \tilde{\mathcal{S}}} \sum_{a \in \mathcal{A}} \mathcal{X}_{ia}^{\omega t} &= \sum_{i \in \tilde{\mathcal{S}}}\theta_i, && \forall \omega \in \Omega \\
       &= 1, && \forall \omega \in \Omega \quad (\textit{Since $\theta$ is a valid PMF over $\tilde{\mathcal{S}}$})
    \end{align*}
\item for $2 \leq t \leq T-1$:
    
    We are given that $\sum\limits_{a \in \mathcal{A}}\mathcal{X}^{\omega t}_{ia}=\sum\limits_{h \in \tilde{\mathcal{S}}}\sum\limits_{a \in \mathcal{A}}\mathcal{X}^{\omega,t-1}_{ha}P_{hai}^\omega$ for each $i \in \tilde{\mathcal{S}}$ and $\omega \in \Omega$ by (\ref{RNEF_flow2}). Moreover, we can re-write $Z^{\omega t}$ as shown in (\ref{RNEF_flow4}). Then, following equations are satisfied for each scenario $\omega \in \Omega$:
    \begin{align*}
        \sum\limits_{i \in \tilde{\mathcal{S}}} \sum\limits_{a \in \mathcal{A}} \mathcal{X}_{ia}^{\omega t} + Z^{\omega t} &= \sum\limits_{i \in \tilde{\mathcal{S}}}\sum\limits_{h \in \tilde{\mathcal{S}}}\sum\limits_{a \in \mathcal{A}}\mathcal{X}^{\omega,t-1}_{ha}P_{hai}^\omega + \sum\limits_{h \in \tilde{\mathcal{S}}}\sum\limits_{a \in \mathcal{A}}\mathcal{X}^{\omega,t-1}_{ha}Q_{ha}^\omega + \mathcal{Z}^{\omega,t-1}\\
         &= \sum\limits_{h \in \tilde{\mathcal{S}}}\sum\limits_{a \in \mathcal{A}}\mathcal{X}_{ha}^{\omega,t-1}\Big[\underbrace{\sum\limits_{i \in \tilde{\mathcal{S}}}P_{hai}^\omega+Q_{ha}^\omega}_{1, \forall (h,a) \in (\tilde{\mathcal{S}},\mathcal{A})}\Big]+\mathcal{Z}^{\omega,t-1}\\
         &= \sum\limits_{h \in \tilde{\mathcal{S}}}\sum\limits_{a \in \mathcal{A}}\mathcal{X}_{ha}^{\omega,t-1}+\mathcal{Z}^{\omega,t-1}
    \end{align*}
     Now, we should prove that $\sum\limits_{h \in \tilde{\mathcal{S}}}\sum\limits_{a \in \mathcal{A}}\mathcal{X}_{ha}^{\omega,t-1}+\mathcal{Z}^{\omega,t-1}$ must be equal to 1. To that end, we utilize \textit{proof by induction} as follows:
     \begin{itemize}
        \item \textit{Base Case: } for $t=2$, $\sum\limits_{h \in \tilde{\mathcal{S}}}\sum\limits_{a \in \mathcal{A}}\mathcal{X}_{ha}^{\omega,t-1}+\mathcal{Z}^{\omega,t-1}\stackrel{?}{=}1, \forall \omega \in \Omega$ 
        
        For $t=2$, the lhs is reduced to $\sum\limits_{h \in \tilde{\mathcal{S}}}\sum\limits_{a \in \mathcal{A}}\mathcal{X}_{ha}^{\omega,1}+\mathcal{Z}^{\omega,1}$. We have already proved that $\sum\limits_{h \in \tilde{\mathcal{S}}}\sum\limits_{a \in \mathcal{A}}\mathcal{X}_{ha}^{\omega,1}=1, \forall \omega \in \Omega$. We also know that $\mathcal{Z}^{\omega,1}=0, \forall \omega \in \Omega$. Thus, base case is satisfied.
        
        \item \textit{Induction Step: } $\sum\limits_{h \in \tilde{\mathcal{S}}}\sum\limits_{a \in \mathcal{A}}\mathcal{X}_{ha}^{\omega,t-1}+\mathcal{Z}^{\omega,t-1}\stackrel{?}{=}1$ given that $\sum\limits_{h \in \tilde{\mathcal{S}}}\sum\limits_{a \in \mathcal{A}}\mathcal{X}_{ha}^{\omega,t-2}+\mathcal{Z}^{\omega,t-2}=1, \forall \omega \in \Omega$
        
        We are given that $\sum\limits_{a \in \mathcal{A}}\mathcal{X}_{ha}^{\omega,t-1}=\sum\limits_{i \in \tilde{\mathcal{S}}}\sum\limits_{a \in \mathcal{A}}\mathcal{X}^{\omega,t-2}_{ia}P_{iah}^\omega, \forall \omega \in \Omega$ for each $h \in \tilde{\mathcal{S}}$ by (\ref{RNEF_flow2}).  Moreover, we can re-write $Z^{\omega,t-1}$ as expressed in (\ref{RNEF_flow4}). Then following equations are satisfied for each scenario $\omega \in \Omega$,
        \begin{align*}
        \sum\limits_{h \in \tilde{\mathcal{S}}} \sum\limits_{a \in \mathcal{A}} \mathcal{X}_{ha}^{\omega,t-1} + \mathcal{Z}^{\omega,t-1} &= \sum\limits_{i \in \tilde{\mathcal{S}}}\sum\limits_{a \in \mathcal{A}}\mathcal{X}_{ia}^{\omega,t-2}\Big[\underbrace{\sum\limits_{h \in \tilde{\mathcal{S}}}P_{iah}^\omega+Q_{ia}^\omega}_{1, \forall (i,a) \in (\tilde{\mathcal{S}},\mathcal{A})}\Big]+\mathcal{Z}^{\omega,t-2}\\
         &= \sum\limits_{i \in \tilde{\mathcal{S}}}\sum\limits_{a \in \mathcal{A}}\mathcal{X}_{ia}^{\omega,t-2}+\mathcal{Z}^{\omega,t-2}\\
         &=1
    \end{align*}
    \end{itemize}

    \item for $t=T:$
    
    We are given that $\mathcal{Y}_i^{\omega} = \sum\limits_{h \in \tilde{\mathcal{S}}}\sum\limits_{a \in \mathcal{A}}\mathcal{X}^{\omega,T-1}_{ha}P_{hai}^\omega$ for each $i \in \tilde{\mathcal{S}}$ and $\omega \in \Omega$. Moreover, we can re-write $Z^{\omega,T-1}$ as depicted in (\ref{RNEF_flow4}). Then, following equations are satisfied for each scenario $\omega \in \Omega$:
    \begin{align*}
        \sum_{i \in \tilde{\mathcal{S}}}\mathcal{Y}_{i}^\omega + \mathcal{Z}^{\omega,T} &= \sum_{i \in \tilde{\mathcal{S}}} \sum\limits_{h \in \tilde{\mathcal{S}}}\sum\limits_{a \in \mathcal{A}}\mathcal{X}^{\omega,T-1}_{ha}P_{hai}^\omega + \sum\limits_{h \in \tilde{\mathcal{S}}}\sum\limits_{a \in \mathcal{A}}\mathcal{X}^{\omega,T-1}_{ha}Q_{ha}^\omega + \mathcal{Z}^{\omega,T-1} \\
          &= \sum\limits_{h \in \tilde{\mathcal{S}}}\sum\limits_{a \in \mathcal{A}}\mathcal{X}_{ha}^{\omega,T-1}\Big[\underbrace{\sum\limits_{h \in \tilde{\mathcal{S}}}P_{hai}^\omega+Q_{ha}^\omega}_{1, \forall (h,a) \in (\tilde{\mathcal{S}},\mathcal{A})}\Big]+\mathcal{Z}^{\omega,T-1} \\
          &= \sum\limits_{h \in \tilde{\mathcal{S}}}\sum\limits_{a \in \mathcal{A}}\mathcal{X}_{ha}^{\omega,T-1}+\mathcal{Z}^{\omega,T-1}\\
          &= 1
    \end{align*}
\end{itemize}

\subsection{Proposition \ref{VI:2}}

\begin{align*}
    \sum_{i\in \tilde{\mathcal{S}}}\mathcal{X}_{i,1}^{\omega t} &\leq \frac{C_t}{n} & \forall t \in \tilde{\mathcal{T}}, \omega \in \Omega\\
    \sum_{t\in \tilde{\mathcal{T}}}\sum_{i\in \tilde{\mathcal{S}}}\mathcal{X}_{i,1}^{\omega t}&\leq \frac{\sum\limits_{t\in \tilde{\mathcal{T}}}C_t}{n} & \forall \omega \in \Omega\\
    \sum_{t\in \tilde{\mathcal{T}}}\bigg[\underbrace{\sum_{i\in \tilde{\mathcal{S}}}\Big[\mathcal{X}_{i,1}^{\omega t}+\mathcal{X}_{i,0}^{\omega t}\Big]+\mathcal{Z}^{\omega t}}_{1 \text{ by } (\ref{VI:1})}\bigg]&\leq \frac{\sum\limits_{t\in \tilde{\mathcal{T}}}C_t}{n} + \sum\limits_{t \in \tilde{\mathcal{T}}} \Big [\sum\limits_{i \in \tilde{\mathcal{S}}}\mathcal{X}^{\omega t}_{i,0} + \mathcal{Z}^{\omega t}\Big]  & \forall \omega \in \Omega\\
    T-1 &\leq \frac{\sum\limits_{t\in \tilde{\mathcal{T}}}C_t}{n} + \sum\limits_{t \in \tilde{\mathcal{T}}} \Big [\sum\limits_{i \in \tilde{\mathcal{S}}}\mathcal{X}^{\omega t}_{i,0} + \mathcal{Z}^{\omega t}\Big]  & \forall \omega \in \Omega\\
     T - \frac{\sum\limits_{t\in \tilde{\mathcal{T}}}C_{t}+n}{n} & \leq \sum\limits_{t\in \tilde{\mathcal{T}}}\sum\limits_{i\in \tilde{\mathcal{S}}}\mathcal{X}^{\omega t}_{i,0}+\mathcal{Z}^{\omega t}  & \forall \omega \in \Omega
\end{align*}

\section{Hierarchical Rules for Transition Probabilities and Rewards}
\label{appendix:montecarlo_transitionprobs}

\subsection{Transition Probabilities}

\begin{enumerate}
    \item Consider two patients in the same awareness level. The patient with worse health status is more likely to worsen.
    \begin{align*}
        P_{1,a,2} &\geq P_{0,a,1} && \forall a \in \{0,1\}\\
        P_{4,a,5} &\geq P_{3,a,4} && \forall a \in \{0,1\}
    \end{align*}
    \item Consider two patients in the same complexity level. The patient with the higher awareness level is less likely to worsen.
    \begin{align*}
        P_{0,a,1} &\geq P_{3,a,4} && \forall a \in \{0,1\}\\
        P_{1,a,2} &\geq P_{4,a,5} && \forall a \in \{0,1\}
    \end{align*}
    \item Consider two patients in the same complexity and awareness levels. The patient taking normal care is more likely to worsen than the patient taking special care.
    \begin{align*}
        P_{0,0,1} &\geq P_{0,1,1}\\
        P_{1,0,2} &\geq P_{1,1,2}\\
        P_{3,0,4} &\geq P_{3,1,4}\\
        P_{4,0,5} &\geq P_{4,1,5}
    \end{align*}
    \item Probability of making a transition from low awareness level to high awareness level are equal regardless of the health status for all patients under same action. Moreover, special care is increasing this probability.
    \begin{align*}
        P_{0,a,3}=P_{1,a,4}=P_{2,a,5} && \forall a \in \{0,1\}\\
        P_{i,1,i+3} \geq P_{i,0,i+3} && \forall i \in \{0,1,2\}
    \end{align*}
    \item Probability of making a transition from high awareness level to low awareness level are equal regardless of the health status for all patients under same action. Moreover, special care is decreasing this probability.
    \begin{align*}
        P_{3,a,0}=P_{4,a,1}=P_{5,a,2} && \forall a \in \{0,1\}\\
        P_{i,0,i-3} \geq P_{i,1,i-3} && \forall i \in \{3,4,5\}
    \end{align*}
    \item A patient cannot change its health status and awareness level in one decision epoch. That is, diagonal arcs in the network is not possible for all patients under any course of action. Therefore, related transition probabilities equal to zero.
    \begin{align*}
        P_{0,0,4}=P_{1,0,5}=P_{3,0,1}=P_{4,0,2}=P_{4,0,0}=P_{5,0,1}=P_{1,0,3}=P_{2,0,4}=0 \\
        P_{0,1,4}=P_{1,1,5}=P_{3,1,1}=P_{4,1,2}=P_{4,1,0}=P_{5,1,1}=P_{1,1,3}=P_{2,1,4}=0
    \end{align*}
    \item There is not any direct transition from \textit{simple} to \textit{complex}.
    \begin{align*}
        P_{0,a,2}=P_{3,a,5}=0 && \forall a \in \{0,1\}
    \end{align*}
    \item Recovery is not possible for all patients under any course of action. Therefore, transition probabilities associated with the arcs to the left equal to zero.
    \begin{align*}
        P_{1,a,0}=P_{2,a,1}=P_{4,a,3}=P_{5,a,4}=0 && \forall a \in \{0,1\}
    \end{align*}
    \item Following relations related to the death probabilities have to be satisfied:
    \begin{align*}
        Q_{2,a}\geq Q_{5,a} = Q_{1,a}\geq Q_{4,a} = Q_{0,a} \geq Q_{3,a} && \forall a \in \{0,1\}
    \end{align*}
    \item Special care reduces the death probabilities
    \begin{align*}
        Q_{i,1}\leq Q_{i,0} && \forall i \in \{0,\dotsc, 5\}
    \end{align*}
    \item There is an upper bound for death probabilities
    \begin{align*}
        Q_{ia} \leq 0.20 && \forall i \in \{0,\dotsc, 5\}, a \in \mathcal{A}
    \end{align*}
\end{enumerate}

\subsection{Rewards}

\begin{enumerate}
    \item Following relations have to be satisfied:
    \begin{align*}
        r_{2,a} \geq r_{5,a}=r_{1,a} \geq r_{4,a} = r_{0,a} \geq r_{3,a} && \forall a \in \{0,1\}
    \end{align*}
    \item Special care is more desirable than normal care.
    \begin{align*}
        r_{i,1} \geq r_{i,0} && \forall i \in \{0, \dotsc, 5\}
    \end{align*}
    \item Final stage rewards are computed as follows:
    \begin{align*}
        R_i \triangleq \frac{\sum\limits_{a \in \mathcal{A}}r_{ia}}{|\mathcal{A}|} && \forall i \in \{0, \dotsc, 5\}
    \end{align*}
    \item Domain of rewards is $[100,1000]$ except $R^D$ where $R^D=0$.
    \begin{align*}
        100 \leq r_{ia} \leq 1000 && \forall i \in \{0, \dotsc, 5\}, a \in \mathcal{A} 
    \end{align*}
\end{enumerate}

\section{Pseudocode of the Strategy Evaluation Algorithm}
\label{appendix:mrp}

\begin{algorithm}[ht!]
\scriptsize
\caption{ Strategy Evaluation Algorithm}
\begin{algorithmic}[1]
\Require $\hat{\Pi}$: the given strategy
\For{$t \in \mathcal{T}$}
    \If{$t=1$}
        \For{$j \in \tilde{\mathcal{S}}$}
            \If{$\Hat{\pi}^{1}_j=1$}
            \State $\mathcal{X}_{j,0}^{\omega,1}\gets0$ \textbf{ and } $\mathcal{X}_{j,1}^{\omega,1}\gets\theta_j\quad \forall \omega \in \Omega$
            \Else
            \State
            $\mathcal{X}_{j,0}^{\omega,1}\gets\theta_j$ \textbf{ and } $\mathcal{X}_{j,1}^{\omega,1}\gets0 \quad \forall \omega \in \Omega$
            \EndIf
            \State $\mathcal{Z}^{\omega,1} \gets 0 \quad \forall \omega \in \Omega$
        \EndFor
    \EndIf
    \If{$1<t<T$}
        \For{$j \in \tilde{\mathcal{S}}$}
        \If{$\Hat{\pi}^{t}_j=1$}
        \State $\mathcal{X}_{j,0}^{\omega t}\gets0$ 
        \State
        $\mathcal{X}_{j,1}^{\omega t}\gets\sum\limits_{i \in \tilde{\mathcal{S}}}\mathcal{X}_{i \pi_i^{t-1}}^{\omega,t-1}P_{iaj}^{\omega}\quad \forall \omega \in \Omega$
        \Else
        \State
        $\mathcal{X}_{j,1}^{\omega t}\gets0$
        \State
        $\mathcal{X}_{j,0}^{\omega t}\gets\sum\limits_{i \in \tilde{\mathcal{S}}}\mathcal{X}_{i \pi_i^{t-1}}^{\omega,t-1}P_{iaj}^{\omega}\quad \forall \omega \in \Omega$
        \EndIf
        \State $\mathcal{Z}^{\omega t}\gets\mathcal{Z}^{\omega, t-1}+ \sum\limits_{i \in \tilde{\mathcal{S}}}\mathcal{X}_{i \pi_i^{t-1}}^{\omega,t-1}Q_{ia}^{\omega}\quad \forall \omega \in \Omega$
        \EndFor
    \EndIf
    \If{$t=T$}
    \For{$j \in \tilde{\mathcal{S}}$}
    \State $\mathcal{Y}^{\omega}_j \gets \sum\limits_{i \in \tilde{\mathcal{S}}}\mathcal{X}_{i \pi_i^{t-1}}^{\omega,t-1}P_{iaj}^{\omega}\quad \forall \omega \in \Omega$
    \EndFor
    \EndIf
\EndFor
\Ensure $(\mathcal{X},\mathcal{Y}, \mathcal{Z})^{\Hat{\Pi}}$: the resulting occupancy measures
\end{algorithmic}\label{Algorithm:MRP}
\end{algorithm}

\section{Evaluating $f^\omega(\Omega)$}
\label{appendix:evaluate_f}

An optimal strategy for a single scenario does not guarantee the feasibility under a large number of scenarios. Because the occupancy measures implied by the strategy of a single scenario may violate the capacity constraints for one of the other scenarios. Thus, we embrace the following approach to evaluate $f^\omega (\Omega)$.

Let $\bar{\Pi}^\omega \triangleq \{\bar{\pi}_i^{\omega t}\}_{i \in \tilde{\mathcal{S}}, t \in \tilde{\mathcal{T}}}$ be the optimal strategy, and $\bar{\pi}_i^{\omega t}$ be the optimal policy at 
stage $t \in \tilde{\mathcal{T}}$ for state $i \in \tilde{\mathcal{S}}$ under scenario $\omega$. We first evaluate $\Bar{\Pi}^\omega$ under $\Omega$ by using Algorithm \ref{Algorithm:MRP}. If it does not involve any infeasibility, the total expected reward characterized by $\bar{\Pi}^\omega$ becomes $f^\omega (\Omega)$. Otherwise, $\bar{\Pi}^\omega$ consists of inapplicable policies. Thus, we find the nearest feasible strategy $\Pi^*$ to $\bar{\Pi}^\omega$, and evaluate $f^\omega(\Omega)$ as the total expected reward characterized by $\Pi^*$. Here, the nearest feasible strategy refers to the one that entails feasible policies with a minimum number of changes in the original policy. The distance is specified based on the changes in the policy since we assume that policy makers tend to apply as few as possible changes in the original plan as inconveniences appear during the execution. In this regard, the nearest strategy $\Pi^*$ to $\bar{\Pi}^\omega$ is determined by the following mixed-integer program. 

\begin{align}
            & {\textbf{min}}
            & &  \sum_{t \in \tilde{\mathcal{T}}} \sum_{i \in \tilde{\mathcal{S}}} \beta_i^t\\
            & \textbf{st;}
            & & (\ref{RNEF_flow1})-(\ref{RNEF_domain4})\\
            &&& \pi_i^t - \bar{\pi}_i^{\omega t} \leq \beta_i^t & \forall t \in \tilde{\mathcal{T}},  i \in \tilde{\mathcal{S}} \\
             &&& \bar{\pi}_i^{\omega t}-\pi_i^t  \leq \beta_i^t & \forall t \in \tilde{\mathcal{T}},  i \in \tilde{\mathcal{S}} \\
            &&& \beta_{i}^{t} \in \{0,1\} & \forall t \in  \tilde{\mathcal{T}}, i \in  \tilde{\mathcal{S}}
\end{align}

\end{document}